\documentclass[a4paper,11pt]{amsart} 
\usepackage{amssymb}
\usepackage{amsmath}
\usepackage{amsthm}
\usepackage{amsbsy}
\usepackage{bbm}
\usepackage{enumerate}
\usepackage{bm}
\usepackage{hyperref}
\date{\today}
\usepackage{cite}
\usepackage{array}
\usepackage{xcolor}
\usepackage{tikz}
\usepackage{graphics}
\usepackage{subcaption}
\usepackage{multirow}
\usepackage{mathtools}
\usepackage{tikz}
\usepackage{caption}
\usepackage{subcaption}
\usepackage{arydshln}
\usepackage{soul}
\usepackage{ulem}
\usetikzlibrary{shapes,snakes}
 \theoremstyle{plain} 
    \newtheorem{theorem}{Theorem}[]
    \newtheorem{lemma}[theorem]{Lemma}
    \newtheorem{proposition}[theorem]{Proposition}
    \newtheorem{corollary}[theorem]{Corollary}

    \newtheorem{definition}[theorem]{Definition}
    
    \newtheorem{result}[theorem]{Result}
    \newtheorem{remark}[theorem]{Remark}
    \newtheorem{example}[theorem]{Example}
    \newtheorem{exercise}[theorem]{Exercise}

\numberwithin{equation}{section}

\def \CC{{\mathcal C}}

\def\E{\mathbb{E}}

\def\Z{\mathbb{Z}}

\def\N{\mathbb{N}}

\def\R{\mathbb{R}}

\def\WW{{\mathcal W}}

\def\Z{\mathbb{Z}}

\def\sgn{{\mb{sgn}}}

\def\<{\langle}
\def\>{\rangle}

\def\mb{\mbox}

\def\bar{\overline}

\newcommand\mnote[1]{} 
\newcommand\be{\begin{equation*}}

\newcommand\ee{\end{equation*}}

\newcommand\ben{\begin{equation}}
\newcommand\een{\end{equation}}
\newcommand\bes{\begin{eqnarray*}}
\newcommand\ees{\end{eqnarray*}}

\newcommand\bex{\begin{exercise}}
\newcommand\eex{\end{exercise}}
\newcommand\beg{\begin{example}}
\newcommand\eeg{\end{example}}
\newcommand\benu{\begin{enumerate}}
\newcommand\eenu{\end{enumerate}}
\newcommand\beit{\begin{itemize}}
\newcommand\eeit{\end{itemize}}
\newcommand\berk{\begin{remark}}
\newcommand\eerk{\end{remark}}
\newcommand\bdefn{\begin{defintion}}
\newcommand\edefn{\end{definition}}
\newcommand\bthm{\begin{theorem}}
\newcommand\ethm{\end{theorem}}
\newcommand\bprf{\begin{proof}}
\newcommand\eprf{\end{proof}}
\newcommand\blem{\begin{lemma}}
\newcommand\elem{\end{lemma}}

\newcommand{\supp}{\mbox{\rm supp}}

\newcommand{\var}{\mbox{\rm Var}}

\newcommand{\Cov}{\mbox{\rm Cov}}
\newcommand{\sm}{{\raise0.3ex\hbox{$\scriptstyle \setminus$}}}

\def\mb{\mbox}

\def\CHI{\mathchoice%
{\raise2pt\hbox{$\chi$}}%
{\raise2pt\hbox{$\chi$}}%
{\raise1.3pt\hbox{$\scriptstyle\chi$}}%
{\raise0.8pt\hbox{$\scriptscriptstyle\chi$}}}
\def\smalloplus{\raise1pt\hbox{$\,\scriptstyle \oplus\;$}}
	\newcommand{\RSC}{Y_d(n,\mathbf{p})}
	
    \newcommand{\domc}{\operatorname{Dom}(c)}
\newcommand{\Ber}[1]{\operatorname{Ber}(#1)}

\title{Semicircle law for multi-parameter random simplicial complexes}
\begin{document}
\author{Kartick Adhikari}
\address{Department of Mathematics\\
 Indian Institute of Science Education and Research Bhopal\\
  Bhauri, Bhopal, Madhya Pradesh 462066, India }
\email{kartick [at] iiserb.ac.in}

\author{Kiran Kumar}
\address{Division of Science\\
        New York University Abu Dhabi\\
         Saadiyat Island, Abu Dhabi, United Arab Emirates}
 \email{kiran.kumar.as.math@gmail.com}

\author{Koushik Saha}
\address{Department of Mathematics\\
	Indian Institute of Technology Bombay\\
	Powai, Mumbai, Maharashtra 400076, India}

\email{koushik.saha [at] iitb.ac.in}

\begin{abstract}



In this paper, we consider the multi-parameter random simplicial complex model, which generalizes the Linial–Meshulam model and random clique complexes by allowing simplices of different dimensions to be included with distinct probabilities. For $n,d \in \N$ and $\mathbf{p}=(p_1,p_2,\ldots, p_d)\in (0,1]^d$, the multi-parameter random simplicial complex $Y_d(n,\mathbf{p})$ is constructed inductively. Starting with $n$ vertices, edges (1-cells) are included independently with probability $p_1$, yielding the Erd\H{o}s-R\'enyi graph $\mathcal{G}(n,p_1)$, which forms the $1$-skeleton. Conditional on the $(k-1)$-skeleton, each possible $k$-cell is included independently with probability $p_k$, for $2 \le k \le d$.
    
We study the signed and unsigned adjacency matrices of $d$-dimensional multi-parameter random simplicial complexes $\RSC,$ under the assumptions $\min_{i=1,\ldots d-1}\liminf p_i >0$ and $np_d \to \infty$ with $p_d=o(1)$. In general, these matrices have random dimensions and exhibit dependency among their entries. We prove that the empirical spectral measures of both matrices converge weakly to the semicircle law in probability. We also prove that the limiting spectral measure of the signed and unsigned adjacency matrices of another associated  model, the multi-parameter upper model, is the semicircle law, under the assumption $np_d(1-p_d) \to \infty$. Further, we  establish convergence results for the number of $(d-1)$-cells and the number of maximal $(d-1)$-cells of multi-parameter simplicial complexes.
\end{abstract}
\maketitle


\section{Introduction and Background}\label{sec:intro}
Random graphs since their introduction in the late 1950s by Erd\H{o}s and R{\'e}nyi \cite{Erdos_Renyi_59,Erdos_Renyi60}, have been widely used in modeling random interactions. Since the advent of this millennium, higher-order interactions have been found to contribute to the behaviors of large real-world systems ranging from the brain\cite{giusti2016two_brain}
to the universe \cite{vandeWeygaert2011}, and simplicial complexes have been successfully used to model them.
The study of random simplicial complexes is motivated both by the need for null models in topological data analysis\cite{Bobrowski_survey} and by the efforts to investigate how the classical properties of random graphs extend to their higher-dimensional analogues.

The first model of random simplicial complexes was introduced by Linial and Meshulam \cite{Linial_RSC_first} for dimension $d=2$, and later this model was generalized to higher dimensions by Meshulam and Wallach \cite{Meshulam_wallach}. In this model, the random simplicial complex has a complete $(d-1)$-skeleton, and $d$-cells are chosen independently with probability $p$. This model is now commonly referred to as the Linial-Meshulam complex. Another important model of random simplicial complexes is the random flag complex (also known as the random clique complex) introduced in \cite{Kahle2009}, which is obtained by considering an Erdős–Rényi random graph and adding all its cliques to the simplicial complex. There also exist several models of geometric random simplicial complexes; we do not go into their details here. For detailed surveys on different models of simplicial complexes, we refer the readers to \cite{Kahle_survey,Bobrowski_survey}.

The multi-parameter random simplicial complex, introduced in \cite{Costa_Farber_multipara1}, allows cells of different dimensions to be included with distinct probabilities. Let $n,d \in \N$ and $\mathbf{p}=(p_1,p_2,\ldots, p_d) \in (0,1]^d$. Given the vector of probabilities $\mathbf{p}$, the multi-parameter random simplicial complex $Y_d(n,\mathbf{p})$ is defined as follows: First, consider $n$ vertices and add edges to them independently with probability $p_1$, i.e., generate an Erd\H{o}s-R{\'e}nyi graph $\mathcal{G}(n, p_1)$ graph. Second, considering this graph as the 1-skeleton of a 2-complex, go over all admissible 2-cells (i.e., 3-cliques) in it, and add them to the complex independently with probability $p_2$. Continue in this fashion; given the $(k-1)$-skeleton, go over all $k$-cells (i.e. $(k+1)$-cliques) such that all of their subsets are contained in the complex and add each of them independently with probability $p_k$ until $k=d$. The resulting simplicial complex is denoted by $\RSC$. For simplicity, we shall use $Y$ to denote $\RSC$, when $n$ and $\mathbf{p}$ are clear from the context.

An alternative construction of the multi-parameter random simplicial complex $\RSC$ is through random hypergraphs. Fix the parameter vector $\mathbf{p}=(p_1,p_2,\ldots, p_d)$ and consider the vertex set $\{1,2,\ldots , n\}$. Generate a random hypergraph, $\mathcal{X}(n ,\mathbf{p})$, by choosing each $j$-cell (in hypergraph terminology, edge of cardinality $j+1$) independently with probability $p_j$, for $1 \leq j \leq d$. Then $\RSC$ has the same distribution as the largest simplicial complex contained in $\mathcal{X}(n ,\mathbf{p})$, and therefore this model is also sometimes called \textit{multi-parameter lower model.} Both the Linial-Meshulam complex and the random flag complex are special cases of the multi-parameter lower model, for the parameter vectors $\mathbf{p}=(1,1,\ldots,1,p)$ and $\mathbf{p}=(p,1,1,\ldots,1)$, respectively. In this paper, by a multi-parameter random simplicial complex, we always mean the multi-parameter lower model.

On a similar vein, the \textit{multi-parameter upper model}, $Z_d(n,\mathbf{p}),$ is defined as the smallest simplicial complex containing $\mathcal{X}(n ,\mathbf{p})$, i.e., 
$Z_d(n,\mathbf{p})$ is the collection of all $\sigma^\prime \subseteq \sigma$ such that $\sigma \in \mathcal{X}(n ,\mathbf{p})$,
\begin{equation}\label{eqn:upper_model_defn}
Z_d(n,\mathbf{p})= \bigcup_{\sigma \in \mathcal{X}(n ,\mathbf{p})}\{ \sigma^\prime: \sigma^\prime \subseteq \sigma\}.  
\end{equation}

The upper model is also the Alexander dual of the multi-parameter lower model \cite[Proposition 4.3]{Farber_Multi_duality}, and their homology groups are related by the famous Alexander duality theorem \cite[Chapter V.4]{Edelsbrunner_book}. Note that here too, the parameter vector $\mathbf{p}=(1,1,\ldots,1,p)$ gives the Linial-Meshulam complex, and therefore this model is another generalization of Linial-Meshulam complex. 




Topological invariants such as Betti numbers, higher-dimensional connectivity and cohomology groups are important tools in topological data analysis, and understanding these properties for random simplicial complexes has significant consequences in topological data analysis. As matrices successfully encode these properties, the study of topological properties often goes hand in hand with the study of random matrices associated with the random simplicial complex. In particular, the local weak limit of the Laplacian matrices has been used in \cite{linial_peled} and \cite{Kahle_flag_complex_Ann.Math} to study the limiting distribution of Betti numbers in Linial-Meshulam complex and random flag complexes, respectively. In \cite{Shu_Kanazawa_Betti}, a similar technique was used to study the limiting distribution of Betti numbers of sparse multi-parameter random simplicial complexes.

The topological properties of the multi-parameter lower model were studied in a series of works by Costa and Farber~\cite{Costa_Farber_multipara1,costafarber_ii,costafarber_iii}.
Most studies on this model have focused on the domain $p_i=n^{-\alpha_i}$ for fixed $\alpha_i \in [0,\infty)$ \cite{Fowler_MultiRSC_homology_2019,kogan2024multidimensional}. In contrast, in this paper, we focus on dense multi-parameter random simplicial complexes ($\liminf p_i>0$ for $1\leq i \leq d-1$ and $np_d(1-p_d) \to \infty$). A closely related domain where $0<\liminf p_i\leq \limsup{p_i}<1 $ for all $i$, is known as the medial domain, and the topological properties of multi-parameter lower model in this domain were recently studied in \cite{Farber_medial1}.
The topological properties of the upper model were studied in \cite{Cooley_Kang_Upper_2020,Farber_Novik_MultiRSC_Upper_2023}. 




In this work, we study the adjacency matrices of multi-parameter random simplicial complexes. We consider two versions of adjacency matrices for simplicial complexes: the signed adjacency matrix and the unsigned adjacency matrix.  For a simplicial complex, its Laplacian matrix is defined as the product of the $d$-boundary operator and its transpose, i.e. $L=\partial_d \, \partial_d^{\, t}$. This generalizes the notion of the Laplace matrix of a graph given by $L=NN^t$, where $N$ is the incidence matrix of any directed graph obtained from $G$
by an arbitrary choice of orientation of the edges \cite[Chapter 1]{Brouwer_Spectra_book}. This definition, among other things, helps in generalizing the matrix-tree theorem to simplicial complexes \cite{Simplicial_matrix_tree}. The signed adjacency matrix is defined using the Laplace matrix as $A^+=D - L$, where $D$ is the degree operator of $(d-1)$-cells (the degree of a
$j$-cell is the number of $(j+1)$-cells that contain it) . The entries of the signed adjacency matrix are $0,1$ and $-1$; the unsigned adjacency matrix is obtained by taking the absolute values of each entry of the signed adjacency matrix.
In this paper, we prove semicircle laws for the signed and unsigned adjacency matrices of both the lower and upper models. We also establish  convergence results for the number of $(d-1)$-cells and the number of maximal $(d-1)$-cells for both models.


\section{Main Results and Proof outlines}
We first introduce some basic terminologies related to simplicial complexes, required to state our main results.
\begin{definition}\label{defn: d-cell}
	Let $V$ be a finite set. A simplicial complex $X$ with vertex set $V$ is a collection $X \subset \mathcal{P}(V)$ such that if $\tau \in X$ and $\sigma \subset \tau$, then $\sigma \in X$. An element of the simplicial complex is called a cell. The dimension of a cell $\sigma \in X$ is defined as $|\sigma|-1$, and a cell of dimension $j$ is called a $j$-cell. The dimension of a non-empty simplicial complex is defined as the maximum of the dimensions of its elements. 
	For a simplicial complex $X$ and integer $j \geq -1$, we denote the set of all $j$-cells in $X$ by $X^j$.
\end{definition}

We follow the convention that  all non-empty simplicial complexes contain the null set and the null set is the only cell of dimension $-1$. For a simplicial complex $X$, the cells of dimension 0 are the vertices.

For every finite set $V$ and $d<|V|$, the set 
\[
\{ \sigma \in \mathcal{P}(V):  |\sigma| \leq d+1 \}
\]
is a simplicial complex, and is called the complete $d$-dimensional simplicial complex on $V$ and denoted by $K_d(V)$ or simply $K_d$. On a related note, we say a simplicial complex has a complete $j$-dimensional skeleton if all subsets of $V$ of dimension $j$ are present in the simplicial complex.

We now proceed to define orientation on a $j$-cell. An orientation is an ordering of the vertices of the $j$-cell, and an oriented cell is represented by a square bracket. Two orientations $[x_0,x_1,\ldots ,x_j]$ and $[y_0,y_1,\ldots , y_j]$ of a $j$-cell are said to be equal if the permutation $g$ given by $g(x_i)=y_i$ is an even permutation. Thus, for every $j \geq 1$, each $j$-cell has exactly two orientations. For an oriented $j$-cell $\sigma$, we use $\overline{\sigma}$ to denote the same $j$-cell with the opposite orientation. As an example, for the $j$-cell $\{  x, y, z\}$, $[x,y,z]=[y,z,x]=[z,x,y]$ and $[x,z,y]=[z,y,x]=[y,x,z]$ are the two different orientations.

Suppose $V$ is an ordered set. The ordering on $V$ induces an ordering on each $j$-cell. The orientation corresponding to this ordering is called the positive orientation, while the other orientation is called the negative orientation. The set of all positively oriented $j$-cells of $X$ is denoted by $X_+^j$ and the set of all oriented $j$-cells is denoted by $X_{\pm}^j$. For $j \geq 1$, we have $|X_{\pm}^j|=2|X^j|$.


\begin{figure}
	\captionsetup[subfigure]{justification=centering}
	\begin{subfigure}[b]{0.45\textwidth}
		\centering
		\includegraphics[height=45mm, width=\textwidth]{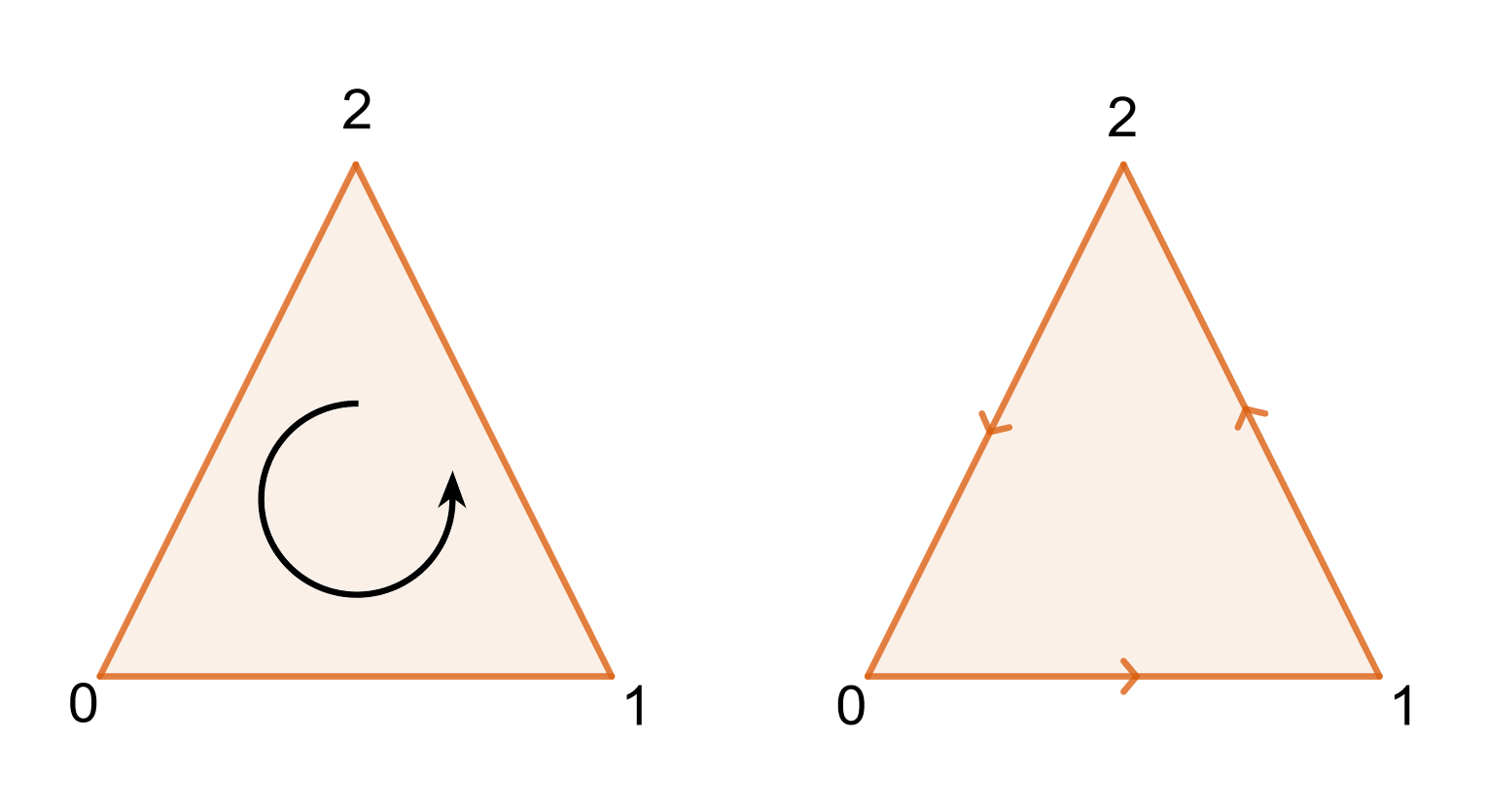}
		\caption{}
	\end{subfigure}
	\hspace{10mm}
	\begin{subfigure}[b]{0.45\textwidth}
		\centering
		\includegraphics[height=45mm, width=\textwidth]{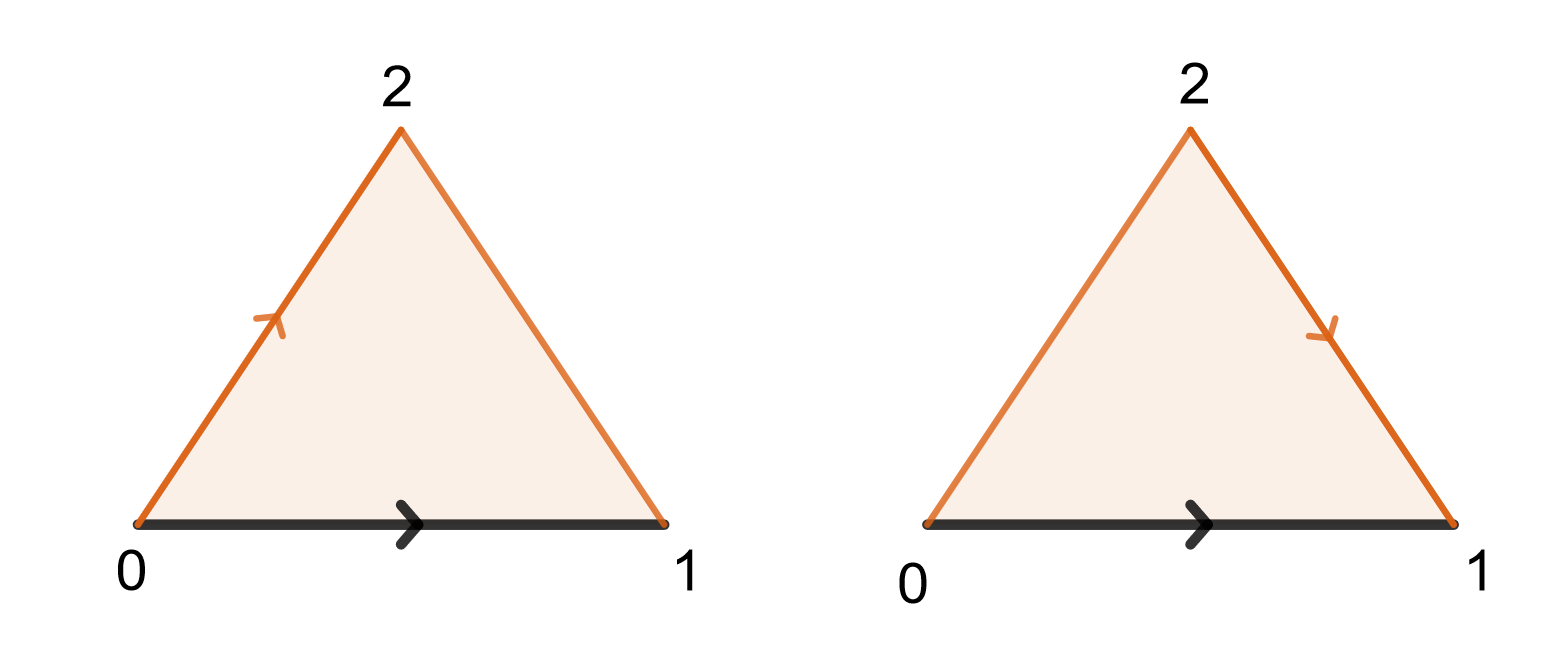}
		\caption{}
	\end{subfigure}
	\caption{Figure (a) shows the orientation induced by the positive orientation of $[0,1,2]$ on its boundary. Figure (b) shows the 1-cell $[0,1]\stackrel{X}{\sim} [0,2]$ and $[0,1]\stackrel{X}{\sim} [2,1]$.}\label{fig:exa orientation}
\end{figure}
For a $j$-cell $\sigma$, we define the boundary of $\sigma$ as $\partial \sigma= \{ \omega \subset \sigma : \text{dim}(\omega)=j-1\}$. An oriented $j$-cell $\sigma=\left[\sigma^{0}, \ldots, \sigma^{j}\right]$ induces an orientation on $\partial \sigma$ as follows: the cell $\left\{\sigma^{0}, \ldots, \sigma^{i-1}, \sigma^{i+1}, \ldots, \sigma^{j}\right\}$ is oriented as $(-1)^i\left[\sigma^{0}, \ldots, \sigma^{i-1}, \sigma^{i+1}, \ldots, \sigma^{j}\right]$, where $-\omega=\bar{\omega}$.


We use the following neighboring relation for oriented cells introduced in \cite{ori_Spectrum_homology}. For a random simplicial complex $X$ of dimension $d$ and $\sigma, \sigma^\prime \in K_{\pm}^{d-1}$, we denote $\sigma \stackrel{X}{\sim} \sigma^{\prime}$ if there exists an oriented $d$-cell $\tau \in X_{\pm}$ such that both $\sigma$ and $\bar{\sigma^{\prime}}$ are in the boundary of $\tau$ as oriented cells. See Figure \ref{fig:exa orientation} (b) for an example.

\begin{definition}\label{defn:adjacency matrix}
	For a simplicial complex $X$ of dimension $d \geq 2$, define
	\begin{enumerate}[(i)]
		\item Unsigned adjacency matrix ($A$): As the $X^{d-1} \times X^{d-1}$ matrix defined by
		\begin{equation*}
			A_{\sigma, \sigma^\prime}= 
			\begin{cases}
				&1 \quad\text{if} \,\,\sigma \cup \sigma^\prime \in X^d, \\
				&0 \quad \text{otherwise.}
			\end{cases}
		\end{equation*}
		\item Signed adjacency matrix ($A^+$): As the $X^{d-1}_+ \times X^{d-1}_+$ matrix defined by
		\begin{equation*}
			A^+_{\sigma, \sigma^\prime}= 
			\begin{cases}
				&1 \quad\text{if} \,\,\sigma \stackrel{X}{\sim} \sigma^\prime , \\
				-&1 \quad \text{if}\,\, \sigma \stackrel{X}{\sim} \bar{\sigma^\prime} , \\
				&0 \quad \text{otherwise.}
			\end{cases}
		\end{equation*}
	\end{enumerate}
\end{definition}
The definition used here is rather simplistic. Nevertheless, it serves our purpose as we are interested in studying the adjacency matrix as a random matrix. For a more detailed explanation, see \cite{Adhikari_Kiran_Saha_RSC_2023}. 

Observe that for all $j \geq 1$, there exists a natural bijection between $X^{d-1}$ and $X_+^{d-1}$. Therefore, the matrix $A^+$ can be indexed by elements of $X^{d-1}$. Hereafter, we consider $A^+$ as a matrix indexed by elements of $X^{d-1}$.

\subsection{Main Theorems}

In this work, we consider sequences of multi-parameter random simplicial complexes $(\RSC)_n$ and  $(Z_d(n,\mathbf{p}))_n$, as the number of vertices, $n$, grows to infinity. 

Consider a sequence of multi-parameter random simplicial complexes $(\RSC)_n$. Let $A_n$ and $A_n^+$ denote the unsigned and signed adjacency matrices of $\RSC$, respectively. Note that if $p_i<1$ for some $1 \leq i \leq d-1$, the set $Y^{d-1}$(and $Y_+^{d-1}$) is random, and therefore $A_n$ (and $A_n^+$) has a random dimension. We look at the normalized adjacency matrices
\begin{align}\label{eqn:An_normalization}
   \frac{1}{\sqrt{nd \prod_{r=1}^{d} p_r^{\beta_r}(1-p_d)}} A_n \quad\text{and} \quad
   \frac{1}{\sqrt{nd \prod_{r=1}^{d} p_r^{\beta_r}(1-p_d)}} A_n^+ \, ,\quad \text{where } \beta_r =\binom{d}{r} ,
\end{align}
and assuming that the terms in the denominator are non-zero. Also, note that here $p_i$ depends on $n$, but we suppress this dependence in the notation. 
To make the denominator in \eqref{eqn:An_normalization} more compact, we introduce the following notation 
\begin{equation}
       m_n(\mathbf{p}) =nd \prod_{r=1}^{d} p_r^{\beta_r}(1-p_d)\label{eqn:adjaceny_normal}.
\end{equation}

For an $n \times n$ symmetric matrix $M_n$ with eigenvalues $\lambda_1,\lambda_2,\ldots, \lambda_n$, the empirical spectral distribution of $M_n$ is defined as
\begin{equation*}
	F_{M_n}(x)=\frac{1}{n} \sum_{i=1}^{n} \mathbf{1}(\lambda_i \leq x).
\end{equation*} 
The probability measure on the real line corresponding to $F_{M_n}$ is known as the empirical spectral measure of $M_n$, and is given by
\begin{equation}\label{defn:spec_dist}
	\mu_{M_n}=\frac{1}{n} \sum_{i=1}^{n} \delta_{\lambda_i}.
\end{equation}
Note that when $M_n$ is a random matrix, $\mu_{M_n}$ is random, and the expectation of $\mu_{M_n}$ is called the expected empirical spectral measure.
The weak limit of $\mu_{M_n}$ is called the limiting spectral measure of $M_n$. 

 Arguably, the most well-known limiting spectral measure in random matrix theory is the Wigner's semicircle law, which emerges as the limiting spectral measure of the Wigner matrix, the symmetric random matrix with i.i.d. upper triangular entries of mean zero and unit variance. The (standard) semicircle law is given by
 \begin{equation}
     \mu_{sc}(a,b)= \frac{1}{2\pi}\int_a^b \sqrt{4-s^2} \cdot \mathbf{1}_{[-2,2]}(s) ds.
 \end{equation}
 
It follows from Wigner's semicircle law and a rank-inequality argument that when $np(1-p) \to \infty$, the limiting spectral measure of the (normalized) adjacency matrix of the Erd\H{o}s-R{\'e}nyi graph $\mathcal{G}(n,p)$ is also the semicircle law \cite{Bollobas_book_2001}. 
The matrix structure of the adjacency matrices of Linial-Meshulam complexes is different from that of Erd\H{o}s-R{\'e}nyi graphs. In particular, the adjacency matrices of Linial-Meshulam complexes are highly sparse, with each row having at most $(n-d)d$ non-zero entries among its ${n \choose d}$ entries.
In \cite{antti_ron}, Knowles and Rosenthal established the semicircle law for Linial-Meshulam complexes: 

\begin{result}[Theorem 2.5,\cite{antti_ron}]\label{res:LSD_linialMeshulam}
For fixed $d \geq 2$, let $A_n^+$ be the signed adjacency matrix of the Linial-Meshulam complex with parameter $p$, and suppose $\lim_{n \to \infty} np(1-p)=\infty$. Then 
$$\mu_{\frac{1}{\sqrt{ndp(1-p)}}A_n^+} \stackrel{w}{\longrightarrow} \mu_{sc} \text{ almost surely,}$$
where $\mu_{sc}$ is the semicircle law.
\end{result}

In \cite{antti_ron}, the authors only considered the signed adjacency matrices. However, it
can be established from Result \ref{res:LSD_linialMeshulam} that the limiting spectral measure of the unsigned adjacency matrix is also the semicircle law. 

Recall that a Linial-Meshulam complex with parameter $p$ is a multi-parameter lower model with the parameter vector $\mathbf{p}=(1,1,\ldots, p)$ and for in this case $m_n(\mathbf{p})=ndp(1-p)$, the normalization for $A_n^+$ in Result \ref{res:LSD_linialMeshulam}. The major contribution of this paper is a generalization of Result \ref{res:LSD_linialMeshulam} for dense multi-parameter random simplicial complexes when $p_d=o(1)$.

\begin{theorem}\label{thm:semicircle}
	For fixed $d \geq 2$, suppose $\min_{i=1,\ldots,d-1}\liminf p_i>0$,  $np_d(1-p_d) \rightarrow \infty$ and $p_d=o(1)$ and  let $A_n$ and $A_n^+$ denote the unsigned and signed adjacency matrices of $Y_d(n,\mathbf{p})$, respectively. Then  both
    \[
    \mu_{\frac{A_n}{\sqrt{m_n(\mathbf{p})}}} \stackrel{w}{\longrightarrow} \mu_{sc} \quad \text{and} \quad \mu_{\frac{A_n^+}{\sqrt{m_n(\mathbf{p})}}} \stackrel{w}{\longrightarrow}  \mu_{sc}
    \]
   in probability.
\end{theorem}

Recall that the multi-parameter upper model is another generalization of the Linial-Meshulam complex. We prove semicircle law for the upper model as well, providing a direct generalization of Result \ref{res:LSD_linialMeshulam}.
\begin{theorem}\label{thm:semicircle_upper}
For fixed $d \geq 2$, suppose $np_d(1-p_d) \rightarrow \infty$ and let $A_n$ and $A_n^+$ denote the unsigned and signed adjacency matrices of $Z_d(n,\mathbf{p})$, respectively. Then  both
    \[
    \mu_{\frac{A_n}{\sqrt{ndp_d(1-p_d)}}} \stackrel{w}{\longrightarrow} \mu_{sc} \quad \text{and} \quad \mu_{\frac{A_n^+}{\sqrt{ndp_d(1-p_d)}}} \stackrel{w}{\longrightarrow}  \mu_{sc}
    \]
   almost surely.
\end{theorem}

Before going to the other result, we first note the following corollary of Theorem \ref{thm:semicircle} and Theorem \ref{thm:semicircle_upper}.
A $(d-1)$-cell of a simplicial complex $X$ is called maximal if it is not contained in any $d$-cell. We denote the number of maximal $(d-1)$-cells of  $X$ by $N_{d-1}(X)$.
Now observe that, since the rows indexed by maximal $(d-1)$-cells in $A_n$ are zero vectors,  the number of maximal $(d-1)$ cells is bounded by the nullity of $A_n$. Therefore we obtain the following. 
\begin{corollary}
\begin{enumerate}[(i)]
    \item Suppose $\min_{i=1,\ldots,d-1}\liminf p_i>0$, $np_d(1-p_d) \rightarrow \infty$ and $p_d=o(1)$, then for the multi-parameter random simplicial complex $Y_d(n,\mathbf p)$,
    $$\frac{ N_{d-1}(Y) }{| Y^{d-1}|} \rightarrow 0 \text{ in probability,}$$
   where $| Y^{d-1}|$ denotes the number of $(d-1)$ cells  in $Y^{d-1}$.
   \item Suppose $np_d(1-p_d) \rightarrow \infty$, then for the multi-parameter upper model $Z_d(n,\mathbf{p})$,
    $$\frac{ N_{d-1}(Z) }{| Z^{d-1}|} \rightarrow 0 \text{ almost surely.}$$
\end{enumerate}
\end{corollary}

By restricting the value of $p_d$ to a slightly smaller domain, we prove that $N_{d-1}(Y)$ converges almost surely to zero. For functions $f,g: \N \to \R$, we say $f \gg g$ if $\frac{g(n)}{f(n)} \to 0$, as $n \to \infty$.

\begin{theorem}\label{thm:maximal_faces}
\begin{enumerate}[(i)]
    \item Suppose $\min_{i=1,\ldots,d-1}\liminf p_i>0$ and $p_d(n) \gg \log n /n$, then 
 $\RSC$ almost surely does not contain any maximal $(d-1)$-cells, i.e., $N_{d-1}(Y) \rightarrow0$ almost surely.
 \item Suppose $p_d(n) \gg \log n /n$, then 
 $Z_d(n,\mathbf{p})$ almost surely does not contain any maximal $(d-1)$-cells, i.e., $N_{d-1}(Z) \rightarrow0$ almost surely.
\end{enumerate}
\end{theorem}
\subsection{Technical difficulties and proof strategy.}
The study of spectral properties of adjacency matrices of random simplicial complexes is relatively new. The first work in this domain was \cite{antti_ron}, where Knowles and Rosenthal established semi-circle law for Linial-Meshulam complex in the domain $np(1-p) \to \infty$. The adjacency matrices of a general multi-parameter random simplicial complexes are more challenging to study compared to the adjacency matrices of Linial-Meshulam complexes, for the following reason:
 \begin{enumerate}[(i)]
     \item For a general multi-parameter simplicial complex, the dimension of its adjacency matrix, which is the number of $(d-1)$-cells of $\RSC$, is random. This is unlike the usual models encountered in random matrix theory, where the dimension is non-random. The randomness of the dimension affects both the definition of the
    empirical spectral distribution and the associated moment
    calculations, in particular the \(k\)-th moment of ESD given by
    \[
        m_k(\mu_{A_n})
        =\frac{1}{\dim(A_n)}\operatorname{Tr}(A_n^k)
        =\frac{1}{\dim(A_n)}
        \sum_{\substack{
            \sigma_0,\ldots,\sigma_{k-1}
            \in Y^{d-1}} \atop \sigma_k=\sigma_0}
        \prod_{j=0}^{k-1}
        A_n(\sigma_j,\sigma_{j+1}).
    \]
    Thus, unlike in the Erd\H{o}s--R\'enyi and Linial--Meshulam
    models, both the normalization and the set over which
    the trace expansion is performed are random. Therefore, the usual classification and enumeration of closed
    walks, words, or simplicial paths cannot be applied directly.
    
    \item  For the Erd\H{o}s--R\'enyi and Linial--Meshulam models, the higher empirical moments of the uncentered adjacency matrix, under the normalization relevant to the semicircle law, diverge for \(k\geq 3\). The standard proof first establishes the semicircle law for
    $A_n-\mathbb{E}[A_n]$ by the moment method and then proves the result to \(A_n\) using a rank inequality
    \cite{Bollobas_book_2001,antti_ron}. However, in the multi-parameter model, \(\mathbb{E}[A_n]\) is not canonically defined because \(A_n\) has a random dimension. We deal with this by extending the adjacency matrix to a non-random dimension, and taking a random centering. 
    \item For Linial-Meshulam complex, by construction, the $d$-cells occur independently of each other.  As a result, its adjacency matrix can be written as
    \[
        A_n=\sum_{\tau\in K^d}\xi_\tau H_\tau,
    \]
    where the matrices \(H_\tau\) are deterministic and the coefficients \(\xi_\tau\) are independent Bernoulli variables.
    For the multi-parameter lower model, since the lower-dimensional skeleton is random, the indicators of distinct \(d\)-cells are generally dependent, with the covariance between entries of the unsigned adjacency matrix given by
\begin{align*}
&\quad \operatorname{Cov}(A_{\sigma_1^{(1)} \sigma_2^{(1)}},A_{\sigma_1^{(2)} \sigma_2^{(2)}})\\
&=\begin{cases}
    (\prod_{\tau \subset S_1}\E X_{\tau})\left( \operatorname{Var}\left(\prod_{\tau \subset S_2} X_{\tau}\right)\right) &\text{if } \left(\sigma_1^{(1)} \cup \sigma_2^{(1)} \right) \cap \left(\sigma_1^{(2)} \cup \sigma_2^{(2)}\right) \neq \phi\\
    0 &\text{otherwise},
\end{cases}
\end{align*}
where $S_1=\mathcal P\left(\sigma_1^{(1)} \cup \sigma_2^{(1)} \right) \Delta \mathcal P\left(\sigma_1^{(2)} \cup \sigma_2^{(2)}\right)$ and $S_2=\mathcal P\left(\sigma_1^{(1)} \cup \sigma_2^{(1)} \right) \cap \mathcal P \left(\sigma_1^{(2)} \cup \sigma_2^{(2)}\right)$, $\Delta$ denotes the symmetric difference of sets and $\mathcal P(S)$ denotes the power set of $S$.
The dependence of entries makes the direct moment calculations further harder.
 \end{enumerate}
\vskip3pt
\noindent\textbf{Proof Strategy.}
 Among the two models considered in this paper, establishing the semicircle law for the multi-parameter lower model (Theorem \ref{thm:semicircle}) is more challenging.  
We provide a broad outline of the strategy employed to prove Theorem \ref{thm:semicircle} and its connection to earlier works. 
 Here, we consider the adjacency matrix of the multi-parameter model and extend it to a matrix of dimension $\binom{n}{d} \times \binom{n}{d}$, called the extended adjacency matrix and denoted by $\widehat{A}_n$,  by adding zero entries everywhere else. However, centering by the expectation of the extended adjacency matrix is not useful in this setting, since the moments of the ESD of the centered matrix also diverge. Instead, we express the extended adjacency matrix $\widehat{A}_n$ as a Hadamard product of $d$ random matrices and center only one of them. We denote this partially  centered version of $\widehat{A}_n$ by $B_n$.
We then study the limiting spectral distribution (LSD) of
${B}_n$, and prove that on a domain, denoted by $\domc$ (see Definition \ref{defn:domc}), the ESD of properly scaled $B_n$ converges almost surely to the probability mixture $\mu_{sc}^{(c)}$, defined by
\begin{equation}\label{convolution_semicircle}
\mu_{sc}^{(c)}([a,b])= 
(1-c)\mathbf{1}(0\in[a,b])
+c\mu_{sc}([a,b]).
\end{equation}
The proof relies on the observation that, under appropriate conditioning, the moment computations become similar to those in \cite{antti_ron}. 

The next step is to relate the LSDs of $B_n$ and $\widehat{A}_n$. For both the Erdős–Rényi graph and the Linial-Meshulam complex, this is straightforward as  $\widehat{A}_n-B_n =\E A_n$ is deterministic for both these models. But for the general model we have $\widehat{A}_n-B_n$   is a random matrix.
Here, we use a Frobenius norm argument to show that the LSD of the extended adjacency matrix is same as that of the centered matrix, making crucial use of the fact that $p_d=o(1)$. 

Finally, we show that, on the domain $\domc$, the proportion of $(d-1)$-cells converges in probability to $c$ (Lemma \ref{lem:num_simplices}). Consequently, we prove that the zero entries added to the extended matrix contribute exactly a mass of $(1-c)\delta_0$ to the LSD of the extended matrix. Therefore, it follows that the remaining mass, $c\mu_{sc}([a,b])$, corresponds to the LSD of the adjacency matrix. After an appropriate normalization, we conclude that, on $\domc$, the LSD of the adjacency matrix is  the semicircle law in probability. Then, we complete the proof of Theorem \ref{thm:semicircle} using a subsequence argument.
\vskip3pt

 \noindent\textbf{Outline of the paper.} In Section \ref{sec:ch3_combina_simplcial}, we introduce the combinatorial objects and preliminary lemmas required for the proofs of our main theorems. In Section \ref{sec:LSD_Hn}, we prove that for a sequence $(\RSC)_n$ belonging to the domain $\domc$, the empirical spectral measure of an appropriately normalized version of $B_n$ converges almost surely to the probability measure $\mu_{sc}^{(c)}$. In Section \ref{sec:main proof}, we prove Theorem \ref{thm:semicircle}, and in Section \ref{sec:upper_model}, we prove Theorem \ref{thm:semicircle_upper}.  Finally, in Section \ref{sec:count_simplices}, we prove Lemma \ref{lem:num_simplices} and Theorem \ref{thm:maximal_faces}.  

 

\subsection{Notations and conventions}
\begin{enumerate}[(i)]
    \item $[n]=\{1,2,\ldots, n\}$.
    \item $|S|$: cardinality of the set S.
    \item $d$: dimension of the multi-parameter random simplicial complex.
    \item $\mathbf{p}$: vector of dimension $d$.
    \item $p_i$: $i$-th component of $\mathbf{p}$.
    \item $\beta_r=\binom{d}{r}$.
    \item $m_n(\mathbf{p})=nd \prod_{r=1}^{d} p_r^{\beta_r}(1-p_d)$.
    \item  $\RSC$: multi-parameter lower model of dimension $d$ with parameter $\mathbf{p}$.
    \item $Z_d(n,\mathbf{p})$: multi-parameter upper model of dimension $d$ with parameter $\mathbf{p}$.
    \item $X^{j}$: set of $j$-cells of a simplicial complex $X.$
    \item $X_+$: set of positively oriented cells of a simplicial complex $X$.
    \item $K_d$: the complete simplicial complex (on $n$ vertices) of dimension $d$.
    \item $N_{d-1}(X)$: number of maximal $(d-1)$-cells of $X$.
    \item $A_n, A_n^+$: signed and unsigned adjacency matrices.
    \item $\widehat{A}_n, \widehat{A}^+_n$: extended unsigned and signed adjacency matrices.
    \item $\mu_n \stackrel{w}{\longrightarrow} \mu$: $\mu_n$ converges to $\mu$ weakly.
\end{enumerate}

\section{Combinatorial objects and Preliminary lemmas}\label{sec:ch3_combina_simplcial}

Recall the definition of a cell from Definition \ref{defn: d-cell}. For $j \geq -1$, we denote the set of all $j$- cells on $\{1,2,\ldots,n\}$ by $K^j$. Here, the dependence on $n$ is suppressed for notational convenience.
\begin{definition}\label{defn: word}
	We define an element of $K^{d-1}$ as a letter. A word of length $k \geq 1$ is a sequence $\sigma_{1} \sigma_{2} \ldots \sigma_{k}$ of letters such that $\sigma_{i} \cup \sigma_{i+1}$ is a $d$-cell. For a word $w=\sigma_{1} \sigma_{2} \ldots \sigma_{k}$, we define
	\begin{equation*}
		\operatorname{supp}_u(w)=\{\tau:\tau \subseteq \sigma_i \cup \sigma_{i+1} \text{ for some }1 \leq i \leq k-1 \text{ and } \operatorname{dim}(\tau)=u\}.
	\end{equation*}
	Note that $\operatorname{supp}_{0}(w)=\sigma_{1} \cup \sigma_{2} \cup \cdots \cup \sigma_{k}$ and $\operatorname{supp}_{d}(w)=\left\{\sigma_{i} \cup \sigma_{i+1} : 1 \leq i \leq k-1\right\}.$ The set $\left\{\sigma_{i}, \sigma_{i+1}\right\}$ is called an edge. We denote by $N_{w}(e)$, the number of times an edge $e$ is crossed by $w$. For  $\tau \in K^d$, we define $N_{w}(\tau):=\sum_{e \in E_{w}(\tau)} N_{w}(e)$, where $E_{w}(\tau)=\left\{\left\{\sigma_{i}, \sigma_{i+1}\right\}: \sigma_{i} \cup \sigma_{i+1}=\tau\right\}$.
\end{definition}

Consider the word $w=\{1,2\}\{1,3\}\{3,4\}\{1,4\}\{2,4\}\{1,2\}$. For the word $w$, we have $\operatorname{supp}_0(w)=\{1,2,3,4\}$, $\operatorname{supp}_1(w)=\{\{1,2\},\{1,3\},\{2,3\},\{1,4\},\{3,4\},\{2,4\}\}$ and $\operatorname{supp}_d(w)=\{\{1,2,3\},\{1,3,4\},\{1,2,4\}\}$. Here, note that each edge appears only once and therefore $N_w(e)=1$ for all edges $e$. Now, for $d$-cells $\tau \in \operatorname{supp}_d(w)$, $N_w(\tau)$ denotes the number of $1 \leq i \leq 5$ such that $\sigma_i \cup \sigma_{i+1}=\tau$, and therefore $N_w(\{1,2,3\})=1$, $N_w(\{1,3,4\})=2$, $N_w(\{1,2,4\})=2$.
\begin{definition}\label{defn:equi_word}
	We say two words $w=\sigma_{1} \sigma_{2} \ldots \sigma_{k}$ and $w^\prime=\sigma_{1}^\prime \sigma_{2}^\prime \ldots \sigma_{k}^\prime$ are equivalent if there exists a bijection $\pi:\operatorname{supp}_0(w) \rightarrow \operatorname{supp}_0(w^\prime)$ such that $\pi(\sigma_i)=\sigma_i^\prime$ for all $i$ and $\pi|_{\sigma_1}$ is a strictly increasing function. 
\end{definition}

Consider words $w=\{1,2\}\{1,3\}\{3,4\}\{2,4\}\{1,2\}$ and $w^\prime=\{3,5\}\{1,3\}\{1,7\}\{5,7\}\{3,5\}$. Here $\operatorname{supp}_0(w)=\{1,2,3,4\}$ and $\operatorname{supp}_0(w^\prime)=\{3,5,1,7\}$. The words $w$ and $w^\prime$ are equivalent, and the appropriate bijection $\pi$ is given by $\pi(1)=3,\pi(2)=5,\pi(3)=1,\pi(4)=7$ and this is the only possible bijection showing equivalence between $w$ and $w^\prime$.

\begin{definition}\label{defn: W_s^k}
	A word $ w=\sigma_{1} \sigma_{2} \ldots \sigma_{\ell}$ is said to be closed if $\sigma_{\ell}=\sigma_{1}$. We denote by $\mathcal{W}_s^k(n,d),$ the set of equivalence classes of closed words of length $k+1$ and $|\operatorname{supp}_0(w)|=s$ such that $N_w(\tau) \neq 1$ for every $\tau \in K^d$. When $n,d$ are clear from the context, we denote the set $\mathcal{W}_s^k(n,d)$ by $\mathcal{W}_s^k$.
\end{definition}

Consider a word $w$ such that $|\operatorname{supp}_0(w)|=s$; we use the following convention to denote a representative element of $[w]$. The first $(d-1)$-cell is always taken as $\{1,2,\ldots ,d\}$ and the subsequent new $0$-cells appear in ascending order from $\{d+1,d+2,\ldots ,s\}$. For example, the representative element for the word $\{3,4\}\{3,6\}\{6,5\}\{3,6\}\{3,4\}$ would be $\{1,2\}\{1,3\}\{3,4\}\{1,3\}\{1,2\}$.

We now recall two lemmata from \cite{Adhikari_Kiran_Saha_RSC_2023}, which will be used to prove Theorem \ref{thm:semicircle}.

\begin{lemma}\cite{Adhikari_Kiran_Saha_RSC_2023}\label{lem: supp_0,supp_d}
	For every word $w$, $\left|\operatorname{supp}_{0}(w)\right| \leq\left|\operatorname{supp}_{d}(w)\right|+d$.
\end{lemma}


\begin{lemma}\cite{Adhikari_Kiran_Saha_RSC_2023}\label{claim:size W_k,s}
	For every $d \geq 2$, $k \geq 2$ and $ d+1 \leq s \leq \bigl\lfloor\frac{k}{2}\bigr\rfloor+d$,
	$$|\mathcal{W}_s^k| \leq  (ds)^{k}.$$
\end{lemma}

Before going to the  next lemma, we make a small comment on the intuition behind Lemma \ref{lem: supp_0,supp_d}. Note that for a word $w=\sigma_1\sigma_2\ldots \sigma_{k+1}$,  the number of $0$-cells in $\sigma_1$, is $d$, whereas there are no $d$-cells in $\sigma_1$. Lemma \ref{lem: supp_0,supp_d} says that the addition of each $0$-cell should result in an addition of a $d$-cell, and therefore we cannot add more than one $0$-cell at any time.

The next lemma gives the number of $j$-cells in the words in $\WW_{\frac{k}{2}+d}^k$, and we see that this is the same for all words $w$ in the set.

\begin{lemma}\label{lemma:card_supp}
	For even $k$ and all $w \in \mathcal{W}_{k/2+d}^k$, $|\operatorname{supp}_u(w)|= {d \choose u+1} + \frac{k}{2}{d \choose u}$ for all $1 \leq u \leq d$.
\end{lemma}
\begin{proof}
	Let $w=\sigma_{1}\sigma_{2}\ldots \sigma_{k+1} \in \WW_{k/2+d}^k$. We count the number of $0$-cells and $d$-cells according to the chronological order of their appearance in the word $w$, proceeding from left to right. For $\sigma_1 \cup \sigma_2$, the number of $u$-cells is ${d+1 \choose u+1}$.  By Lemma \ref{lem: supp_0,supp_d}, $\sigma_j \cup \sigma_{j+1}$ is a new $d$-cell if and only if $\sigma_{j+1}$ contains a new $0$-cell. Suppose all the $u$-cells in $\sigma_{1}\cup \sigma_2, \sigma_2 \cup \sigma_{3},\ldots, \sigma_{i-1} \cup \sigma_{i}$ are counted. Suppose the $d$-cell $\sigma_i \cup \sigma_{i+1}$ has already been counted. Then, clearly, all the $u$-cells in $\sigma_i \cup \sigma_{i+1}$ have been counted. Next, suppose the $d$-cell $\sigma_i \cup \sigma_{i+1}$ is appearing for the first time. Then from Lemma \ref{lem: supp_0,supp_d}, it follows that $\sigma_i \cup \sigma_{i+1}$ contains new $0$-cells. Note that some of the $u$-cells in $\sigma_i \cup \sigma_{i+1}$ are already counted, and the ones that are not counted are exactly the ones containing the new $0$-cells. Therefore, the number of new $u$-cells added is ${d \choose u}$. Finally, note that the number of times a new $d$-cell appears is $k/2-1$, and therefore
	\begin{equation*}
		|\operatorname{supp}_u(w)|={d+1 \choose u+1}+\left(\frac{k}{2}-1\right){d \choose u}={d+1 \choose u+1}-{u \choose d}+\frac{k}{2}{d \choose u}={d \choose u+1}+\frac{k}{2}{d \choose u}.
	\end{equation*}
\end{proof}

%

\section{Limiting spectral distribution of centered Hadamard product}\label{sec:LSD_Hn}

Recall from Definition \ref{defn:adjacency matrix} that, in general, the dimensions of the matrices $A_n$ and $A_n^+$ are random, as $Y^d$ is random. We consider the extended unsigned adjacency matrix $\widehat{A}_n$  indexed by $\sigma, \sigma^\prime \in K^{d-1}$ and defined by

\begin{equation}\label{eqn:adj_hadamard_extended}
    \left(\widehat{A}_n\right)_{\sigma \sigma^\prime} =\begin{cases}
        \left(A_n\right)_{\sigma \sigma^\prime}    &\text{ if } \sigma, \sigma^\prime \in Y^{d-1} \\
        0 &\text{ otherwise.}
    \end{cases}
\end{equation}

Since the indexing is by the set of all $(d-1)$-cells, the dimension of the matrix is $\binom{n}{d}$, which is non-random. 

Recall the random hypergraph $\mathcal{X}(n, \mathbf{p})$ introduced in Section \ref{sec:intro}.
Note that for $\sigma, \sigma^\prime \in K^{d-1}$ such that $\sigma \cup \sigma^\prime \in K^d$, $\sigma \cup \sigma^\prime \in Y^d$ if and only if $\tau \in \mathcal{X}(n,\mathbf{p})$ for all $\tau \subseteq \sigma \cup \sigma^\prime$. Since $\tau$ are chosen independently, the entries of the extended adjacency matrix can be given by
\begin{equation}\label{eqn:extended_adj_entry}
  (\widehat{A}_n)_{\sigma \sigma^\prime}=\prod_{\tau \subset \sigma \cup \sigma^\prime} \chi_\tau,  
\end{equation}
where $\chi_{\tau}$ are independent Bernoulli random variables with parameter $p_{\dim(\tau)}$.
Further, note that if either $\sigma \notin Y^{d-1}$ or $\sigma^\prime \notin Y^{d-1}$, then we have $\chi_\tau=0$ for some $\tau$, and therefore we have $(\widehat{A}_n)_{\sigma \sigma^\prime}=0$ in this case. Also, note that even though we don't explicitly include in the notation, $\chi_\tau$ depends on $n$, since $\mathbf{p}$ depends on $n$.

Using \eqref{eqn:extended_adj_entry}, we can express $\widehat{A}_n$ as a Hadamard product, 
\begin{equation}\label{eqn:adj_hadamard}
    \widehat{A}_n =\odot_{j=1}^d A_n(j) \text{ where } A_n(j)_{\sigma \sigma^\prime}=\begin{cases}
    \prod_{\tau \subseteq \sigma \cup \sigma^\prime \atop \dim(\tau)=j} \chi_{\tau} &\text{ for } \sigma \cup \sigma^\prime \in K^d, \\
    0 &\text{ otherwise}.
        
    \end{cases}
\end{equation}

Similarly, the extended signed adjacency matrix of $\RSC$  can be written as the Hadamard product
$$\widehat{A}_n^+=\sgn(K_d) \odot \widehat{A}_n, $$
where $\sgn(K_d)$ is a matrix with non-random entries given by  
$$\sgn(K_d)_{\sigma \sigma^\prime}= 
			\begin{cases}
				&1 \quad\text{if} \,\,\sigma \stackrel{K_d}{\sim} \sigma^\prime , \\
				-&1 \quad \text{if}\,\, \sigma \stackrel{K_d}{\sim} \bar{\sigma^\prime} , \\
				&0 \quad \text{otherwise.}
			\end{cases}$$
for $\sigma,\sigma^\prime \in {K_+^{d-1}}$.

Observe that when considering a sequence of random simplicial complexes,
$(\RSC)_{n \in \N}$, we work on the $\sigma$-algebra generated by the triangular family of independent random variables $\chi_\tau^{(n)}$ indexed by all $j$-cells $\tau$ ($j \leq d$) on the infinite vertex set $\{1, 2, 3, \ldots\}$,
where $\chi_\tau$ is a Bernoulli random variable.

For $d \in \Z_+, n \geq d+1$ and $\mathbf{p}=(p_1,p_2,\ldots,p_d)$, consider the multi-parameter random simplicial complex $Y_d(n,\mathbf{p})$. Define the matrix $B_n$ as
\begin{equation}\label{eqn:centred_adj}
	B_n= \odot_{j=1}^{d-1}A_n(j) \odot (A_n(d)-\E A_n(d)),
\end{equation}
where $A_n(j)$ are as defined in \eqref{eqn:adj_hadamard}, and the matrix $H_n$ as
$$H_n=\frac{1}{\sqrt{m_n(\mathbf p)}} B_n,$$
where $m_n(\mathbf p)$ as defined in \eqref{eqn:adjaceny_normal}. 
Note that here $\E A_n(d)=p_d \mathbb{A}_n(d)$, where $\mathbb{A}_n(d)$ is the unsigned adjacency matrix of the complete $d$-dimensional simplicial complex on $n$ vertices.

In a similar fashion, for the signed adjacency matrix, we define
\begin{equation*}
	B_n^+=\sgn(K_d)\odot B_n= \sgn(K_d) \odot \left(\odot_{j=1}^{d-1}A_n(j) \odot (A_n(d)-p_dA_n(d))\right). 
\end{equation*}
and its normalized version $H_n^+$ as
$$H_n^+=\frac{1}{\sqrt{m_n(\mathbf p)}} B_n^+.$$
\begin{remark}\label{rem:LM}
 For a Linial-Meshulam complex $Y$, the parameter vector is of the form $\mathbf{p}=(1,\ldots,1,p)$, and therefore, the set of $(d-1)$-cells of $Y$ is non-random with $Y^{d-1}=K^{d-1}$. Moreover, for every $\tau$ with $\dim(\tau) \leq d-1$, $\chi_\tau \equiv 1$. As a consequence, for all $j \leq d-1$, $A_n(j)=\mathbb{A}_n$, the unsigned adjacency matrix of $K_d$, and therefore $A_n=\widehat{A}_n=A_n(d)$ for Linial-Meshulam complex. \vskip2pt
\end{remark}

Our main theorem in this section is about the limiting spectral distribution of $H_n$ and $H_n^+$ on a specific domain, which is defined below.
\begin{definition}\label{defn:domc}
    For $c>0$, we say a sequence of random simplicial complexes $\left( \RSC \right)_{n}$ belongs to $\domc$ if $\min_{i=1,\ldots,d-1}\liminf p_i>0$,
    \begin{equation*}
np_d(1-p_d) \rightarrow \infty \text{ and } p_1^{\beta_2} p_2^{\beta_3}\cdots p_{d-1}^{\beta_d} \rightarrow c,
\end{equation*}
 where $\beta_r=\binom{d}{r}$ for $ 1 \leq r \leq d$.
\end{definition}
Before going to the main theorem of this section we recall $\mu_{sc}^{(c)}$ from \eqref{convolution_semicircle} and observe that $\mu_{sc}^{(c)}$ is the distribution of the random variable $Z=XY$ where $X \sim \mu_{sc}$, $Y \sim \mbox{Ber}(c)$ and $X,Y$ are independent. The moments of $\mu_{sc}^{(c)}$ are therefore given by 
\begin{equation}\label{eqn_moment_mu_sc^c}m_k(\mu_{sc}^{(c)})=c\, m_k(\mu_{sc})=\begin{cases}
        c\, \CC_{k/2} &\text{ if } k \text{ is even,}\\
        0 &\text{ otherwise,}
    \end{cases}\end{equation}
    where $\CC_{k/2}$ is the $k/2$-th Catalan number.

In the following theorem, we compute the limiting spectral measures of the centered adjacency matrices, $H_n$ and $H_n^+$ on $\domc$.
\begin{theorem}\label{thm:LSD_B_n}
	For a sequence $(\RSC)_n$ of multi-parameter random simplicial complexes in $\domc,\ c>0$, the limiting spectral measure of both $H_n$ and $H_n^+$ is $\mu_{sc}^{(c)}$ 
    almost surely.
\end{theorem}

We first show that the expected empirical spectral measure of $H_n$ converges weakly to $\mu_{sc}^{(c)}$.

\begin{lemma}\label{lemma:EESD_unsigned}
	For a sequence $(\RSC)_n$ of multi-parameter random simplicial complexes in $\domc, c>0$, the expected empirical spectral measure of $H_n$ converges weakly to $\mu_{sc}^{(c)}$.
\end{lemma}
\begin{proof}
	For $k=0$, we have $\int_{\mathbb{R}} x^k\mu_{H_n}(dx)=1$.
	For $k \geq 1$, note from \eqref{defn:spec_dist} that
	\begin{equation}\label{eqn:k-th moment exp1}
		\mathbb{E}\left[\int_{\mathbb{R}} x^{k} \mu_{H_n}(dx)\right]=\frac{1}{N(m_n(\mathbf{p}))^{k/2}} \sum_{\sigma_{1}, \ldots, \sigma_{k} \in K^{d-1}} \mathbb{E}\left[B_{\sigma_{1} \sigma_{2}} B_{\sigma_{2} \sigma_{3}} \cdots B_{\sigma_{k-1} \sigma_{k}} B_{\sigma_{k} \sigma_{1}}\right],
	\end{equation}
	where $B_{\sigma_i \sigma_{i+1}}$ is the element of the matrix $B_n$ indexed by $(d-1)$ cells $\sigma_i $ and $\sigma_{i+1}$, and $N=\operatorname{dim}(B_n)={n \choose d}$.
	
By the definition of  $B_n$, $B_{\sigma_i \sigma_{i+1}}$ is non-zero only if $\sigma_i \cup \sigma_{i+1}$ is a $d$-cell. As a result, the summation in (\ref{eqn:k-th moment exp1}) can be restricted to the summation over closed words of length $k+1$. Thus, we get
\begin{equation}\label{eqn:k-th moment exp2}
	\mathbb{E}\left[\int_{\mathbb{R}} x^{k} \mu_{H_n}(dx)\right]=\frac{1}{N(m_n(\mathbf{p}))^{k/2}}\sum_{\substack{\text{closed words } w \\ \text{of length } k+1}} \mathbb{E}\left[B_{\sigma_{1} \sigma_{2}} B_{\sigma_{2} \sigma_{3}} \cdots B_{\sigma_{k-1} \sigma_{k}} B_{\sigma_{k} \sigma_{1}}\right],
\end{equation}
	where $w=\sigma_1\sigma_2\ldots \sigma_{k+1}$.

	Fix a word $w=\sigma_1\sigma_2\ldots \sigma_{k+1}$, consider the event 
    \begin{equation}\label{eqn:E_w_definition}
        E_w=\{A_n(j)_{\sigma_{i}\sigma_{i+1}}=1 \ \forall \, 1 \leq i \leq k \text{ and } 1 \leq j \leq d-1\},
    \end{equation} where $A_n(j)$ are the matrices in \eqref{eqn:adj_hadamard}. Note that the value of the summand in \eqref{eqn:k-th moment exp2} is equal to zero in $E_w^c$ and therefore, for each closed word $w$, the summand in \eqref{eqn:k-th moment exp2} can be replaced with
\begin{equation*}
	\mathbb{E}\left[B_{\sigma_{1} \sigma_{2}} B_{\sigma_{2} \sigma_{3}} \cdots B_{\sigma_{k-1} \sigma_{k}} B_{\sigma_{k} \sigma_{1}}|E_w\right]\mathbb{P}(E_w)
\end{equation*}

Note that if the $d$-cells $\sigma_i \cup \sigma_{i+1}$ and $ \sigma_j 
\cup \sigma_{j+1}$ are distinct, the random variables $A_n(d)_{\sigma_i \sigma_{i+1}}$ and $A_n(d)_{\sigma_j \sigma_{j+1}}$ on $E_w$ are i.i.d. with their common distribution $\Ber{p_d}$, and denoted by $(\chi_\tau)_{\tau \in K^d}$. Recalling $N_w(\tau)$ from Definition \ref{defn: word},
 we get that if $N_w(\tau)=1$ for some $\tau \in \operatorname{supp}_d(w)$, then $\mathbb{E}\left[(\chi_\tau-p_d)^{N_{w}(\tau)}\right]=0$. Also, $\mathbb{E}\left[(\chi_\tau-p_d)^{N_{w}(\tau)}\right]=1$ if $N_w(\tau)=0$ for some $\tau \in K^d$. Thus, the summation in (\ref{eqn:k-th moment exp2}) can be restricted to summation over closed words $w$ of length $k+1$ such that $N_w(\tau) \neq 1$ for all $\tau \in K^d$. 
\begin{equation}\label{eqn:k-th moment exp3}
	\mathbb{E}\left[\int_{\mathbb{R}} x^{k} \mu_{H_n}(dx)\right] =\frac{1}{N(m_n(\mathbf{p}))^{k/2}}\sum_{\substack{\text{closed words } w \\ \text{of length } k+1 \\ {N_w(\tau) \neq 1 \forall \tau \in K^d}}} \mathbb{E}\left[B_{\sigma_{1} \sigma_{2}} B_{\sigma_{2} \sigma_{3}} \cdots B_{\sigma_{k} \sigma_{k+1}}|E_w\right] \mathbb{P}(E_w).
\end{equation}	
Consider the Linial-Meshulam complex determined by the random variables $(X_\tau)_{\tau \in K^d}$, and let $M$ denote the centered unsigned adjacency matrix of this Linial-Meshulam complex. Note that on the set $E_w$, $B_{\sigma_{1} \sigma_{2}} B_{\sigma_{2} \sigma_{3}} \cdots  B_{\sigma_{k} \sigma_{k+1}}=M_{\sigma_{1} \sigma_{2}} M_{\sigma_{2} \sigma_{3}} \cdots  M_{\sigma_{k} \sigma_{k+1}}$, and therefore the right hand side of \eqref{eqn:k-th moment exp3} can be expressed as
\begin{equation}\label{eqn:k-th moment exp4}
	\frac{1}{N(m_n(\mathbf{p}))^{k/2}}\sum_{s=d+1}^{\lfloor\frac{k}{2} \rfloor+d}\sum_{w \in \WW_s^k} \mathbb{E}\left[M_{\sigma_{1} \sigma_{2}} M_{\sigma_{2} \sigma_{3}} \cdots M_{\sigma_{k-1} \sigma_{k}} M_{\sigma_{k} \sigma_{k+1}}\right]\mathbb{P}(E_w),
\end{equation}	
where $\WW_s^k$ is as defined in Definition \ref{defn: W_s^k}.  

For a closed word $w$, define
	\begin{equation*}
		T(w):= \mathbb{E}\left[M_{\sigma_{1} \sigma_{2}} M_{\sigma_{2} \sigma_{3}} \cdots M_{\sigma_{k-1} \sigma_{k}} M_{\sigma_{k} \sigma_{k+1}}\right]= \prod_{\tau \in K^{d}} \mathbb{E}\left[(\chi_\tau-p_d)^{N_{w}(\tau)}\right],
	\end{equation*}
where $\{\chi_\tau\}_{\tau \in K^d}$ are i.i.d. Bernoulli random variables with parameter $p_d$.
Notice that $T(w^{\prime})= T(w)$ for all $w^{\prime} \sim w$ and further $\mathbb{P}(E_w)=\mathbb{P}(E_{w^\prime})$. 
	
By the computation in the proof of Lemma 3.2 of \cite{antti_ron}, it follows that the contribution for $s<k/2+d$ in \eqref{eqn:k-th moment exp4} is of the order $o(1)$. Hence, 

\begin{equation*}
\lim_{n \rightarrow \infty} \mathbb{E}\left[\int_{\mathbb{R}} x^{k} \mu_{H_n}(dx)\right]=0 \text{ for odd $k$ and }
\end{equation*}
\begin{equation*}
	\lim_{n \rightarrow \infty} \mathbb{E}\left[\int_{\mathbb{R}} x^{k} \mu_{H_n}(dx)\right] =\lim_{n \rightarrow \infty} \frac{1}{N(m_n(\mathbf{p}))^{k/2}}\sum_{w \in \mathcal{W}_{k/2+d}^k} T(w)|[w]|\mathbb{P}(E_w) \text{ for even $k$,}
\end{equation*}
where $[w]$ denotes the equivalence class of $w$.
Next, we compute $\mathbb{P}(E_w)$. For simplicity, we use $Y$ to denote the random simplicial complex $\RSC$. Note that
\begin{align*}
	&\quad \mathbb{P}(E_w) \\
    &=\mathbb{P}(\operatorname{supp}_{d-1}(w)\subseteq Y^{d-1}|\operatorname{supp}_{d-2}(w)\subseteq Y^{d-2})
			 \times \mathbb{P}(\operatorname{supp}_{d-2}(w)\subseteq Y^{d-2}|\operatorname{supp}_{d-3}(w)\subseteq Y^{d-3}) \\ &\hspace{20mm} \times \cdots \times \mathbb{P}(\operatorname{supp}_{1}(w) \subseteq Y^1) \\
			&=p_{d-1}^{|\operatorname{supp}_{d-1}(w)|} p_{d-2}^{|\operatorname{supp}_{d-2}(w)|} \cdots p_1^{|\operatorname{supp}_1(w)|} \\
			&=\prod_{i=1}^{d-1}p_i^{\beta_{i+1}+\frac{k}{2}\beta_i}.
\end{align*}
Here, the first equality follows since, given $\supp_{u}(w)$, all elements of $\supp_{u+1}(w)$ are chosen independently with probability $p_{u+1}$, and the last equality follows from Lemma \ref{lemma:card_supp}.

Consider $w \in \WW_{k/2+d}^k$. By Lemma~\ref{lem: supp_0,supp_d},
\[
|\supp_d(w)| \geq |\supp_0(w)| - d
= \left(\frac{k}{2}+d\right)-d
= \frac{k}{2}.
\]
On the other hand, since $N_w(\tau)\neq 1$ for every $\tau\in K^d$, it follows that $|\supp_d(w)|\leq \frac{k}{2}$ and hence, $|\supp_d(w)|=\frac{k}{2}$. Consequently, every $\tau\in \supp_d(w)$ must satisfy
$N_w(\tau)=2.$ Therefore, $N_w(\tau)=2$ for every $w\in \WW_{k/2+d}^k$ and every $\tau\in \supp_d(w)$. As a result, we get $T(w)=(p_d(1-p_d))^{k/2}$. Further from Lemma 3.11 of \cite{antti_ron}, it follows that $|\WW_{k/2+d}^k|=\CC_{k/2}d^{k/2}$, where $\CC_r$ is the $r$-th Catalan number, and by Lemma 20 of \cite{Adhikari_Kiran_Saha_RSC_2023}, $|[w]|=\frac{n!}{(n-k/2-d)!d!}$. Combining all the above, we get the limiting even moments as
\begin{align*}
	\lim_{n \rightarrow \infty}	\mathbb{E}\left[\int_{\mathbb{R}} x^{k} \mu_{H_n}(dx)\right] &= \lim_{n \rightarrow \infty} \frac{\frac{n!}{(n-k/2-d)!d!}}{{n \choose d}}\frac{1}{(nd \prod_{r=1}^{d} p_r^{\beta_r}(1-p_d))^{k/2}}\\  &\times \CC_{k/2}d^{k/2}(p_d(1-p_d))^{k/2}\prod_{i=1}^{d-1}p_i^{\beta_{i+1}+\frac{k}{2}\beta_i} \\
	&=  \CC_{k/2}\lim_{n \rightarrow \infty}\prod_{i=1}^{d-1}p_i^{\beta_{i+1}} = \CC_{k/2}c.
\end{align*}
Since $\CC_{k/2}c$ is the $k$-th moment of $\mu_{sc}^{(c)}$  for even values of $k$ (see  \eqref{eqn_moment_mu_sc^c}), the result follows from the method of moments.
\end{proof}

Next, we prove that the moments of the empirical spectral measure of $H_n$ converge to the corresponding moments of $\mu_{sc}^{(c)}$ almost surely. We prove this using the Borel-Cantelli lemma and the following lemma.

\begin{lemma}\label{lemma: sum Var B_n}
	For $d \geq 2$ and $ n \geq d+1$, let $m_k(H_n)$ denote the $k$-th moment of $H_n$. For every sequence $(\RSC)_n$ of multi-parameter random simplicial complexes in $\domc, c>0$,
	$ \mathbb{E} \left[m_k(H_n)- \mathbb{E}m_k(H_n)\right]^2 =O(1/n^2)$ for every positive integer $k$.
\end{lemma}
\begin{proof}
For all positive integers $k$, we have
	\begin{align*}
		m_k(H_n)&=\frac{1}{N(m_n(\mathbf p))^{k/2}} \sum_{\sigma_{1}, \ldots, \sigma_{k} \in K^{d-1}} B_{\sigma_{1} \sigma_{2}} B_{\sigma_{2} \sigma_{3}} \cdots B_{\sigma_{k-1} \sigma_{k}} B_{\sigma_{k} \sigma_{1}}\nonumber\\
		&=\frac{1}{N(nd \prod_{r=1}^{d} p_r^{\beta_r}(1-p_d))^{k/2}} \sum_{w}  \prod_{i=1}^k B_{\sigma_i \sigma_{i+1}}, 
	\end{align*}
where the summation is over closed words $w=\sigma_1 \sigma_2\ldots \sigma_{k+1}$ of length $k+1$ on the vertex set $[n]$. Note that since $\prod_{r=1}^{d-1}p_r^{\beta_{r+1}} \to c>0$, $\liminf p_i>0$ for all $i$ and therefore $\prod_{r=1}^{d-1}p_r^{\beta_r}$ is bounded away from zero. Hence, proving the lemma is equivalent to showing 
	\begin{align}\label{eqn: m_k(B_n) var}
	\frac{1}{N^2\left(np_d(1-p_d)\right)^k}\displaystyle\sum_{w_1,w_2}\mathbb{E}\left(\prod_{i=1}^k B_{\sigma_i^{(1)} \sigma_{i+1}^{(1)}}B_{\sigma_i^{(2)} \sigma_{i+1}^{(2)}}\right)-\mathbb{E}\left(\prod_{i=1}^k B_{\sigma_i^{(1)} \sigma_{i+1}^{(1)}}\right)\mathbb{E}\left(\prod_{i=1}^k B_{\sigma_i^{(2)} \sigma_{i+1}^{(2)}}\right)
	\end{align}
    is of the order $O(\frac{1}{n^2})$,
	where $w_1=\sigma_1^{(1)} \sigma_2^{(1)}\ldots \sigma_{k+1}^{(1)}$ and $ w_2=\sigma_1^{(2)} \sigma_2^{(2)}\ldots \sigma_{k+1}^{(2)}$ are  closed words of length $k+1$ on the vertex set $[n]$. 
	
	For closed words $w_1,w_2$, consider the event
	$$E_{w_1,w_2}=\{A_n(u)_{\sigma_{i}^{(1)}\sigma_{i+1}^{(1)}}=1=A_n(u)_{\sigma_{i}^{(2)}\sigma_{i+1}^{(2)}}, \ \forall \ 1 \leq u \leq d-1, 1 \leq i \leq k\},$$
	where the matrices $A_n(u)$ are as defined in \eqref{eqn:adj_hadamard}.
	Note that the term in \eqref{eqn: m_k(B_n) var} corresponding to closed words $w_1,w_2$ is equal to zero on $E_{w_1,w_2}^c$.
	As a result, each term in the summation in \eqref{eqn: m_k(B_n) var} is equal to
	\begin{align}\label{eqn:var_Bn_2}
		&\mathbb{E}\left(\prod_{i=1}^k B_{\sigma_i^{(1)} \sigma_{i+1}^{(1)}}B_{\sigma_i^{(2)} \sigma_{i+1}^{(2)}}|E_{w_1,w_2}\right)\mathbb{P}(E_{w_1,w_2}) \nonumber\\&\hspace{20mm}-\mathbb{E}\left(\prod_{i=1}^k B_{\sigma_i^{(1)} \sigma_{i+1}^{(1)}}|E_{w_1}\right)\mathbb{E}\left(\prod_{i=1}^k B_{\sigma_i^{(2)} \sigma_{i+1}^{(2)}}|E_{w_2}\right)\mathbb{P}(E_{w_1})\mathbb{P}(E_{w_2}) \nonumber\\
		&= \prod_{\tau \in K^d} \mathbb{E} \chi_{\tau}^{N_{w_1}(\tau)+N_{w_2}(\tau)}\mathbb{P}(E_{w_1,w_2})-  \prod_{\tau \in K^d} \mathbb{E} \chi_{\tau}^{N_{w_1}(\tau)}  \prod_{\tau \in K^d} \mathbb{E} \chi_{\tau}^{N_{w_2}(\tau)}\mathbb{P}(E_{w_1})\mathbb{P}(E_{w_2}), 
	\end{align}
    where the sets $E_{w_1}, E_{w_2}$ are as defined in \eqref{eqn:E_w_definition}, 
	 $\{\chi_\tau\}_{\tau \in K^d}$ are i.i.d. centered $\Ber{p_d}$ random variables and $N_{w_i}(\tau)$ is as defined in Definition \eqref{defn: word}.  
	
	Suppose for a pair of closed words $w_1,w_2$ and $\tau \in K^d$, $N_{w_1}(\tau) + N_{w_2}(\tau)=1$. Since $\E \chi_\tau=0$, it follows that the term in \eqref{eqn:var_Bn_2} is equal to zero for the pair $w_1,w_2$. Further note that if $\supp_1(w_1) \cap \supp_1(w_2) = \phi$, then $\mathbb{P}(E_{w_1,w_2})=\mathbb{P}(E_{w_1})\mathbb{P}(E_{w_2})$, and as a result, the term in \eqref{eqn:var_Bn_2} is equal to zero here as well. Hence, we need to consider only the tuples of the form $(w_1,w_2)$ with $N_{w_1}(\tau)+N_{w_2}(\tau) \neq 1$ for all $\tau \in K^d$ and $\supp_1(w_1) \cap \supp_1(w_2) \neq \phi$.
	
	 We now define an equivalence relation between the tuples of words $(w_1,w_2)$ and $(w_1^\prime,w_2^\prime)$. Two 2-tuples of words $(w_1,w_2)$ and $(w_1^\prime,w_2^\prime)$ are said to be equivalent if there exists a map $\pi: \operatorname{supp}_0(w_1,w_2) \rightarrow  \operatorname{supp}_0(w_1^\prime
	,w_2^\prime) $ such that the restriction $ \pi\big|_{\operatorname{supp}_0(w_1)}$ defines an equivalence between $w_1$ and $w_1^\prime$ (see Definition \ref{defn:equi_word}), and  $\pi\big|_{\operatorname{supp}_0(w_2)}$ is a bijection from $\operatorname{supp}_0(w_2)$ to $\operatorname{supp}_0(w_2^\prime)$. Note that if the tuples of words $(w_1,w_2)$ and $(w_1^\prime,w_2^\prime)$ are equivalent, then the term in \eqref{eqn:var_Bn_2} corresponding to $(w_1,w_2)$ and $(w_1^\prime,w_2^\prime)$ are equal.

	Hence \eqref{eqn: m_k(B_n) var} can be written as
	\begin{align}\label{eqn:variance_two_terms}
		\frac{1}{N^2\left(np_d(1-p_d)\right)^k} \displaystyle\sum_{s=d+1}^{k+2d} &\sum_{(w_1,w_2) \in W_2(k,s,d)}  \Big(\prod_{\tau \in K^d} \mathbb{E} \chi_{\tau}^{N_{w_1}(\tau)+N_{w_2}(\tau)} \mathbb{P}(E_{w_1,w_2}) \nonumber \\&-  \prod_{\tau \in K^d} \mathbb{E} \chi_{\tau}^{N_{w_1}(\tau)}  \prod_{\tau \in K^d} \mathbb{E} \chi_{\tau}^{N_{w_2}(\tau)}\mathbb{P}(E_{w_1})\mathbb{P}(E_{w_2})\Big)  \big|[(w_1,w_2)]\big|, 
	\end{align}
	where $W_2(k,s,d)$ is the set of all equivalence classes of tuples $(w_1,w_2)$ such that $w_1,w_2$ are closed words of length $k+1$, $|\operatorname{supp}_0(w_1,w_2)|=s, \supp_1(w_1)\cap \supp_{1}(w_2) \neq \phi $ and $N_{w_1}(\tau)+N_{w_2}(\tau) \neq 1$ for all $\tau \in K^d$. Here $[(w_1,w_2)]$ denotes the set of all elements in the equivalence class of $(w_1,w_2)$.
	\vskip3pt
	Now, consider the term 
	\begin{equation}\label{eqn:var_order_first_term}
		\frac{1}{N^2\left(np_d(1-p_d)\right)^k} \displaystyle\sum_{s=d+1}^{k+2d} \sum_{(w_1,w_2) \in W_2(k,s,d)}  \prod_{\tau \in K^d} \mathbb{E} \chi_{\tau}^{N_{w_1}(\tau)+N_{w_2}(\tau)} \mathbb{P}(E_{w_1,w_2}) \big|[(w_1,w_2)]\big|.
	\end{equation}
	As $\chi_\tau$ is a {centered} Bernoulli random variable with parameter $p_d$,
	\begin{align*}
		& \prod_{\tau \in K^d} \mathbb{E}\left[\chi_\tau^{N_{w_1}(\tau)+N_{w_2}(\tau)}\right]\\
        &=\prod_{\substack{\tau \in K^d \\
				N_{w_1}(\tau)+N_{w_2}(\tau)\geq 2}} p_d(1-p_d)\left[(1-p_d)^{N_{w_1}(\tau)+N_{w_2}(\tau)-1}+(-p_d)^{N_{w_1}(\tau)+N_{w_2}(\tau)-1}\right] \\
		& \leq \prod_{\substack{\tau \in K^d \\
				N_{w_1}(\tau)+N_{w_2}(\tau) \geq 2}} p_d(1-p_d)=\left(p_d(1-p_d)\right)^{\left|\operatorname{supp}_d(w_1,w_2)\right|}.\\
	\end{align*}
	
	 We make the following two claims:
	\vskip5pt
	\noindent \textbf{Claim 1:} For closed words $w_1,w_2$ with $|\operatorname{supp}_0(w_1,w_2)|=s$, $\big|[(w_1,w_2)]\big|=O(n^s)$.
	\begin{proof}[Proof of Claim 1]
The number of ways to choose $s$ many 0-cells for $\operatorname{supp}_0(w_1,w_2)$ is $n\choose s$. Now for a fixed choice of $\operatorname{supp}_0(w_1,w_2)$, the number of ways to form closed words $w_1,w_2$ is finite and it does not depend on $n$. Therefore $\big|[(w_1,w_2)]\big|=O(n^s)$.
	\end{proof}
    
	\noindent \textbf{Claim 2:} Let $(w_1,w_2) \in W_2(k,s,d)$ such that $N_{w_1}(\tau)+N_{w_2}(\tau) \neq 1$ for all $\tau \in K^d$ and $\operatorname{supp}_1(w_1) \cap \operatorname{supp}_1(w_2) \neq \phi$, then $$|\operatorname{supp}_0(w_1,w_2)| \leq |\operatorname{supp}_d(w_1,w_2)|+2d-2 .$$
	\begin{proof}[Proof of Claim 2]
		We first count the number of 0-cells and $d$-cells in $w_1$, and note from Lemma \ref{lem: supp_0,supp_d} that $|\operatorname{supp}_0(w_1)|-|\operatorname{supp}_d(w_1)| \leq d$. Let $$\ell=\min \{j: \operatorname{supp}_1(\sigma_j^{(2)} \cup \sigma_{j+1}^{(2)}) \cap \operatorname{supp}_1(w_1) \neq \phi\ \}.$$
        Now we count the number of 0-cells and $d$-cells in $\sigma_{\ell}^{(2)} \cup \sigma_{\ell+1}^{(2)}$. Since $\operatorname{supp}_1(\sigma_\ell^{(2)} \cup \sigma_{\ell+1}^{(2)}) \cap \operatorname{supp}_1(w_1)$ contains at least two $0$-cells that are already counted in $w_1$, the number of new $0$-cells in $\sigma_{\ell}^{(2)} \cup \sigma_{\ell+1}^{(2)}$ is at most $d-1$. Therefore, the difference between the number of $0$-cells and the number of $d$-cells counted in $w_1$ and $\sigma_{\ell}^{(2)} \cup \sigma_{\ell+1}^{(2)}$  is at most $d+(d-2)$.

        To count the remaining $0$-cells and $d$-cells in $w_2$, we start by counting the new $0$-cells in $\sigma_{\ell-1}^{(2)}$ and proceed in descending order till $\sigma_{1}^{(2)}$. Note that if $\tau^\prime=\sigma_{j}^{(2)} \cup \sigma_{j-1}^{(2)}$ contains a new $0$-cell that is not yet counted, then the $d$-cell $\tau^\prime$ has also not been counted yet. After reaching $\sigma_{1}^{(2)}$, we start counting the new $0$-cells in $\sigma_{\ell+2}^{(2)}$ and proceed in ascending order till $\sigma_{k+1}^{(2)}$. Note that in this case, if  $\tau^\prime=\sigma_{j}^{(2)} \cup \sigma_{j+1}^{(2)}$ contains a new $0$-cell, then the $d$-cell $\tau^\prime$ is also new. In short, for each new $0$-cell  in $\operatorname{supp}_0(w_2)\setminus \operatorname{supp}_0(\sigma_\ell^{2}\cup\sigma_{\ell+1}^{(2)})$ which has not appeared in $\operatorname{supp}_0(w_1)$,  a new $d$-cell appears in $\operatorname{supp}_d(w_2)$ associated to it. 
        Therefore,  $|\operatorname{supp}_0(w_1,w_2)| \leq | \operatorname{supp}_d(w_1,w_2)|+d+(d-2)$. 
	\end{proof} 
	
	Note that by Claim 1, $\big|[(w_1,w_2)]\big|=O(n^s)$, and by Claim 2, $| \operatorname{supp}_0(w_1,w_2)| \leq | \operatorname{supp}_d(w_1,w_2)|+2d-2$. 
    Therefore, we get that \eqref{eqn:var_order_first_term} is of the order
	\begin{align*}
		&O\left(\frac{1}{n^{2d}}\right) \times  O\left(\frac{1}{\left(np_d(1-p_d)\right)^k}\right) \times O(n^s) \times O\left(\left(p_d(1-p_d)\right)^{| \operatorname{supp}_d(w_1,w_2)|}\right)\\
        &\leq O\left(\left(n\right)^{s-2d-k}\right)\times O\left(\left(p_d(1-p_d)\right)^{| \operatorname{supp}_0(w_1,w_2)|-2d+2-k}\right)\\
        &=O(n^{-2})\times O\left(\left(np_d(1-p_d)\right)^{s-2d+2-k}\right)\\
        &\leq O(n^{-2}).
	\end{align*}
	The last inequality follows from the observation that $s- k+2d-2\leq 0$, since $w_1$ and $w_2$ are closed words of length $k+1$ and $N_{w_1}(\tau)+N_{w_2}(\tau) \neq 1$ for all $\tau$, and $np_d(1-p_d)\to\infty$ as $n\to\infty$. A similar argument shows that the second term of \eqref{eqn:variance_two_terms} is also of the order $O(n^{-2})$, and this completes the proof of the lemma.
\end{proof}
As a consequence of Lemmata \ref{lemma:EESD_unsigned} and \ref{lemma: sum Var B_n} and Borel-Cantelli lemma, it follows that the limiting spectral measure of $H_n$ is $\mu_{sc}^{(c)}$ almost surely.
\vskip5pt
\noindent\textbf{Proof for the signed version.} 
The proof of the signed case is similar to that of the unsigned case, and here we only sketch the proof. 

First, we prove the convergence of the expected empirical spectral measure of $H_n^+$. Here, we proceed in the same fashion as Lemma \ref{lemma:EESD_unsigned}, and carrying out analogous calculations, we get the equivalent form of \eqref{eqn:k-th moment exp4} as
\[
\frac{1}{N(m_n(\mathbf p))^{k/2}}\sum_{s=d+1}^{\lfloor\frac{k}{2} \rfloor+d}\sum_{w \in \WW_s^k} \widetilde{T}(w) \mathbb{P}(E_w),
\]
where $\displaystyle{\widetilde{T}(w)=  \prod_{\tau \in K^{d}} \mathbb{E}\left[(\chi_\tau-p_d)^{N_{w}(\tau)}\right]\sgn(w,\tau)}$ with $\displaystyle{\sgn(w,\tau)=(-1)^{\sum_{e \in \widehat{E}_w(\tau)}N_w(\tau)}}$ and $\widehat{E}_w(\tau) = \{\, \{\sigma, \sigma'\} \in E_{w}(\tau) : \sigma \overset{K_d}{\sim} \overline{\sigma'} \,\}$, taking care of the sign. Since $|\sgn(w,\tau)|=1$, it follows from the same argument in Lemma \ref{lemma:EESD_unsigned} that the contribution of the terms $w \in \WW_s^k$, for $s<\frac{k}{2}+d$ is of the order $o(1)$, and therefore we have \begin{equation*}
\lim_{n \rightarrow \infty} \mathbb{E}\left[\int_{\mathbb{R}} x^{k} \mu_{H_n^+}(dx)\right]=0 \text{ for odd $k$ and }
\end{equation*}
\begin{equation*}
	\lim_{n \rightarrow \infty} \mathbb{E}\left[\int_{\mathbb{R}} x^{k} \mu_{H_n^+}(dx)\right] =\lim_{n \rightarrow \infty} \frac{1}{N(m_n(\mathbf p))^{k/2}}\sum_{w \in \mathcal{W}_{k/2+d}^k} \widetilde{T}(w)|[w]|\mathbb{P}(E_w) \text{ for even $k$}.
\end{equation*}
For $w \in \WW_k^{\frac{k}{2}+d}$, it was proved in \cite[Lemma 3.11]{antti_ron} that $\sgn(w,\tau)=1$, and we have $T(w)=\widetilde{T}(w)$ for all $w \in \WW_k^{\frac{k}{2}+d}$. Therefore, the even moments of $\mu_{H_n^+}$ also converge to the even moments of $\mu_{sc}^{(c)}$.

To prove the equivalent version of Lemma \ref{lemma: sum Var B_n}, note that for the signed case, the equivalent expression of \eqref{eqn:variance_two_terms} would have sign terms that are either $+1$ or $-1$. As this change does not affect the order of \eqref{eqn:variance_two_terms}, we have that the variance of $m_k(H_n^+)$ is of the order $O(n^{-2})$. This proves that the limiting spectral measure of $H_n^+$ is $\mu_{sc}^{(c)}$ almost surely. This completes the proof of Theorem \ref{thm:LSD_B_n}.

\section{Proof of semicircle law for $p_d=o(1)$}\label{sec:main proof}
The objective of this section is to prove Theorem \ref{thm:semicircle}. 
The proof is organized into the following two steps:

\textbf{Step 1.} For a sequence $(\RSC_n \subset \domc$ with $p_d = o(1)$, we prove that the empirical spectral distribution (ESD) of the normalized extended adjacency matrix $\frac{\widehat{A}_n}{\sqrt{m_n(\mathbf{p})}}$ converges almost surely to $\mu_{sc}^{(c)}$.

\textbf{Step 2.} As a consequence of Step 1, we deduce through a subsequence argument that in the regime $p_d = o(1)$, the ESD of $\frac{A_n}{\sqrt{m_n(\mathbf{p})}}$ converges to the semicircle law in probability.



Before going to the individual steps, we expand the Hadamard product in \eqref{eqn:centred_adj} to write $\widehat{A}_n$ in terms of $B_n$. Recall that the last term in the definition of $B_n$ is $(A_n(d)-\E A_n(d))$ and further note that $\E A_n(d)=p_d\mathbb{A}_n(d)$ where $\mathbb{A}_n(d)$ is the unsigned adjacency matrix of the complete $d$-dimensional simplicial complex on $n$ vertices. Therefore, we have
\begin{align*}
	B_n&= \odot_{j=1}^{d-1}A_n(j) \odot (A_n(d)-p_d\mathbb{A}_n(d))\\
	&= \odot_{j=1}^{d} A_n(j)-\left(\odot_{j=1}^{d-1}A_n(j)\odot p_d \mathbb{A}_n(d)\right)\\
	&= \odot_{j=1}^{d} A_n(j)-p_d\odot_{j=1}^{d-1}A_n(j) \\
	&= \widehat{A}_n - p_d\odot_{j=1}^{d-1}A_n(j).
\end{align*}
That is,
\begin{equation}\label{eqn:An,Bn_relation}
	\widehat{A}_n= B_n +  p_d\odot_{j=1}^{d-1}A_n(j).
\end{equation}
Similarly, for the signed version, we have 
\begin{align*}
	B_n^+&= \sgn(K_d) \odot \left( \odot_{j=1}^{d-1}A_n(j) \odot (A_n(d)-p_d\mathbb{A}_n(d))\right)\\
	&= \sgn(K_d) \odot \left( \widehat{A}_n - p_d\odot_{j=1}^{d-1}A_n(j)\right)\\
    &= \widehat{A}^+_n - p_d\sgn(K_d) \odot\left(\odot_{j=1}^{d-1}A_n(j)\right),
\end{align*}
implying
\begin{equation}\label{eqn:An+,Bn+_relation}
	\widehat{A}^+_n= B_n^+ +  p_d\sgn(K_d) \odot\left(\odot_{j=1}^{d-1}A_n(j)\right).
\end{equation}

\begin{figure}
	\includegraphics[height=100mm, width=150mm]{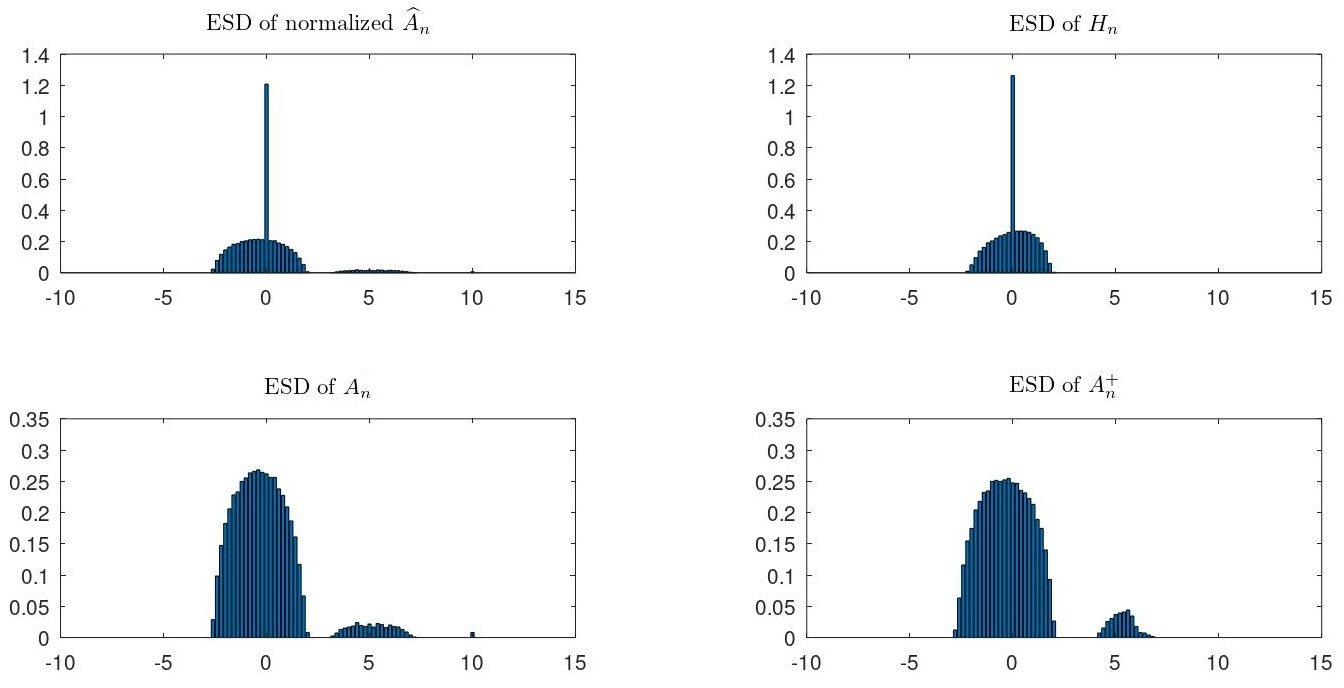}
	\caption{The empirical spectral distributions of: (top-left) normalized extended adjacency matrix $\widehat{A}_n/{\sqrt{\mathbf{m}_n(\mathbf{p})}}$, (top-right) its centered version $H_n$, (bottom-left) the normalized unsigned adjacency matrix $A_n/\sqrt{\mathbf{m}_n(\mathbf{p})}$ and (bottom-right) the normalized signed adjacency matrix $A_n^+/\sqrt{\mathbf{m}_n(\mathbf{p})}$. The histograms are for $n=40,d=2$ and $\mathbf{p}=(0.8,0.7)$, averaged over 10 realizations.}\label{fig:LSD_Bn,An}
\end{figure}
\subsection{Step 1:}\label{sec:LSD_An}
In this step, we consider sequences of random simplicial complexes $\RSC$ belonging to the domain $\domc$ with the additional assumption that $p_d=o(1)$. The main result of this subsection is the following:
\begin{proposition}\label{thm:high_decay_LSD}
	For a sequence $\RSC$ in  $\domc$ such that $p_d=o(1)$, the limiting spectral measure of both $\frac{\widehat{A}_n}{\sqrt{m_n(\mathbf{p})}}$ and $\frac{\widehat{A}^+_n}{\sqrt{m_n(\mathbf{p})}}$ is $\mu_{sc}^{(c)}$  almost surely.
\end{proposition}
	We prove Proposition \ref{thm:high_decay_LSD} using the following lemma which is a generalization of \cite[Lemma 2.4.3]{tao_RMTbook_12}, and follows an analogous proof. 
	\begin{lemma}\label{lemma:Frobenius_bound}
		For any $N \times N$ Hermitian matrices $A$ and $B$, any $\lambda$ and any $\epsilon >0$, we have
		$$
		\mu_{\frac{1}{\sqrt{np_d}}(A+B)}(-\infty, \lambda) \leq \mu_{\frac{1}{\sqrt{np_d}}(A)}(-\infty, \lambda+\epsilon\sqrt{\frac{N}{np_d}})+\frac{1}{\epsilon^2 N^2}\|B\|_F^2
		$$
		and similarly
		$$
		\mu_{\frac{1}{\sqrt{np_d}}(A+B)}(-\infty, \lambda) \geq \mu_{\frac{1}{\sqrt{np_d}}(A)}(-\infty, \lambda-\epsilon\sqrt{\frac{N}{np_d}})-\frac{1}{\epsilon^2 N^2}\|B\|_F^2,
		$$
		where $\|B \|_F$ denotes the Frobenius norm of $B$.
	\end{lemma}
\begin{proof}[Proof of Proposition \ref{thm:high_decay_LSD}]
	We prove the theorem only for the unsigned version, and the proof of the signed version follows along similar lines. 
    
    Define $C_n=p_d \odot_{j=1}^{d-1}A_n(j)$, and choose $\epsilon= \frac{p_d^{3/4}}{n^{\frac{d-1}{2}}}$ and $N=\binom{n}{d}$. 
	Note that $$\epsilon \sqrt{\frac{N}{np_d}}=O\left(\frac{p_d^{1/4}}{n^{\frac{d-1}{2}}} \times n^{d/2} \times n^{-1/2}  \right)= O(p_d^{1/4})=o(1), \mbox{ as }p_d=o(1).$$
	Since the number of non-zero entries of $A_n$ in each row is uniformly bounded above by $nd$ with each entry bounded by 1, we have
	\begin{align*}
		\frac{1}{\epsilon^2 N^2} \|C_n\|_F^2 &\leq \frac{1}{\epsilon^2 N^2} \times N \times n\times d \times p_d^2 \\
		&=O\left( \frac{n^{d-1}}{p_d^{3/2}}\times n^{-d} \times n \times p_d^2 \right) \\
		&=O(\sqrt{p_d})=o(1).
	\end{align*}
	It follows that $\frac{1}{\epsilon^2 N^2} \|C_n\|_F^2$ goes to zero as $n \rightarrow \infty$.
	
    Note that for $(\RSC)_n$ in $\domc$, $\prod_{i=1}^{d-1} p_i^{\beta_i}$ is bounded away from zero, and therefore the conclusion of Lemma \ref{lemma:Frobenius_bound} stays the same even after multiplying the denominator by $\sqrt{\prod_{i=1}^{d-1} p_i^{\beta_i}(1-p_d)}$. 
    Therefore, using Theorem \ref{thm:LSD_B_n} and Lemma \ref{lemma:Frobenius_bound}, we get that for all real $\lambda \neq 0$,
    \begin{multline*}
        \mu_{sc}^{(c)}(-\infty,\lambda) \leq \liminf_n \mu_{\frac{1}{\sqrt{m_n(\mathbf{p})}}(B_n+C_n)}(-\infty,\lambda) \\ \leq\limsup_n \mu_{\frac{1}{\sqrt{m_n(\mathbf{p})}}(B_n+C_n)}(-\infty,\lambda) \leq  \mu_{sc}^{(c)}(-\infty,\lambda) 
    \end{multline*}
	 almost surely. As a result, $ \mu_{\frac{1}{\sqrt{m_n}}(B_n+C_n)}(-\infty,\lambda)$ converges to $ \mu_{sc}^{(c)}(-\infty,\lambda)$ for $\lambda \neq 0$, almost surely. This completes the proof of the proposition. 
\end{proof}

\subsection{Step 2.} In this subsection, we consider random simplicial complexes $(\RSC)_n$ in a larger domain, where $\min_{i=1,\ldots,d-1}\liminf_{n} p_i>0$, $np_d(1-p_d)\to\infty$ and $p_d=o(1)$. We prove that the empirical spectral distribution of the adjacency matrix, rather than the extended adjacency matrix considered in Proposition \ref{thm:high_decay_LSD}, converges to the semicircle law in probability and thus completing  the proof of Theorem \ref{thm:semicircle}.

First we state two lemmas which will be required in the proof of Theorem \ref{thm:semicircle}. 

\begin{lemma}\label{lem:num_simplices}
Suppose $(\RSC)_n$ belongs to $\domc$ for some $c>0$, then $\frac{| Y^{d-1}|}{{n \choose d}} \rightarrow c$ in probability.
\end{lemma}
We prove Lemma \ref{lem:num_simplices} in Section \ref{sec:count_simplices}.
\begin{lemma}\label{lem:subsequence_conv}\cite[Lemma 4.2]{kallenberg2002foundationsofmodernprobability}
Let $X,X_1,X_2,\ldots$ be random elements in a metric space $(S,\rho)$.  
Then $X_n \to X$ in probability if and only if every subsequence $N' \subset \mathbb{N}$ has a
further subsequence $N'' \subset N'$ such that $X_n \to X$ a.s. along $N''$.
\end{lemma}

It is well-known that the space of probability measures is a metric space and therefore by Lemma \ref{lem:subsequence_conv}, to prove that $\mu_{A_n}$ converges weakly to the semicircle law in probability, it is sufficient to show that for any subsequence $(n_k)$, there exists a further subsequence $(n_{k_\ell})$ such that $\mu_{A_{n_{k_\ell}}}$ converges weakly to the semicircle law almost surely, and this is exactly what we will establish.

\begin{proof}[Proof of Theorem \ref{thm:semicircle} ]
Let $(\mathbf{p}_n=(p_1,p_2,\ldots,p_d))_{n\geq1}$ be a sequence of $d$-tuples such that $\min_{i=1,\ldots, d-1} \liminf p_i>0$ and $p_d=o(1)$. Then there exists a subsequence of $\mathbf{p}_n$ belonging to $\domc$, and further by Lemma \ref{lem:num_simplices} and Lemma \ref{lem:subsequence_conv}, there exists a further subsequence $\mathbf{p}_{n_k}$, say,  such that $\frac{| Y_d^{d-1}(n_k,\mathbf p_{n_k})|}{\binom{n_k}{d}}$ converges to $c$ almost surely. Fix this subsequence $\mathbf{p}_{n_k}$.

By Proposition \ref{thm:high_decay_LSD}, for the subsequence $\mathbf{p}_{n_k}$, $\mu_{\frac{\widehat{A}_{n_k}}{\sqrt{m_{n_k}(\mathbf p_{n_k})}}}$ converges weakly to $\mu_{sc}^{(c)}$ almost surely.

Now observe that for each realization of $\RSC$, the matrix $\widehat{A}_n(\omega)$ is a block matrix (with appropriate row and column reordering) of the form
\begin{equation}\label{eqn:block_form}
  \widehat{A}_n(\omega)=\left[
\begin{array}{c:c}
A_n(\omega) & \mathbf{0} \\ \hdashline
\mathbf{0} & \mathbf{0}
\end{array}
\right],  
\end{equation}
where $A_n$ is the unsigned adjacency matrix, which we are interested in. This shows that 
\begin{align*}
F_{\frac{\widehat{A}_n}{\sqrt {m_n(\mathbf p)}}}(x) &=\frac{\operatorname{dim}(A_n)}{\binom{n}{d}}F_{\frac{A_n}{\sqrt{m_n(\mathbf p)}}}(x) \text{ for } x <0  \text{ and}\\
 F_{\frac{\widehat{A}_n}{\sqrt{m_n(\mathbf p)}}}(x) &=\frac{\operatorname{dim}(A_n)}{\binom{n}{d}}F_{\frac{A_n}{\sqrt{m_n(\mathbf p)}}}(x)+ 1- \frac{\operatorname{dim}(A_n)}{\binom{n}{d}}   \text{ for } x \geq 0,
\end{align*}
where $ F_{\widehat{A}_n},F_{A_n}$ denotes the distribution functions of $\mu_{\widehat{A}_n}$ and $\mu_{A_n}$, respectively.

Note that by Proposition \ref{thm:high_decay_LSD}, 
\begin{align*}
&\lim_{k \to \infty} F_{\frac{\widehat{A}_{n_k}}{\sqrt{m_{n_k}(\mathbf p_{n_k})}}}(x) =F_{ \mu_{sc}^{(c)}}(x)=cF_{\mu_{sc}}(x), \ 
\mbox{ almost surely for }x<0,\\  
&\lim_{k \to \infty} F_{\frac{\widehat{A}_{n_k}}{\sqrt{m_{n_k}(\mathbf p_{n_k})}}}(x) =F_{\mu_{sc}^{(c)}}(x)=cF_{\mu_{sc}}(x)+(1-c),\ \mbox{ almost surely for }x > 0,\end{align*} 
where $F_{\mu_{sc}^{(c)}}$ and $F_{\mu_{sc}}$ denote the distribution functions of $\mu_{sc}^{(c)}$ and $\mu_{sc}$ respectively.
On the other hand, by Lemma \ref{lem:num_simplices}, $\frac{\dim(A_{n_k})}{\binom{n_k}{d}} \to c$ almost surely, and therefore \begin{align*}
&\lim_{k\to\infty}F_{\frac{\widehat{A}_{n_k}}{\sqrt{m_{n_k}(\mathbf p_{n_k})}}}(x) = c \lim_{k \to \infty } F_{\frac{A_{n_k}}{\sqrt{m_{n_k}(\mathbf p_{n_k})}}}(x), \ \mbox{ almost surely for }x<0,\\
& \lim_{k \to \infty}F_{\frac{\widehat{A}_{n_k}}{\sqrt{m_{n_k}(\mathbf p_{n_k})}}}(x) = 
c \lim_{n \to \infty } F_{\frac{A_{n_k}}{\sqrt{m_{n_k}(\mathbf p_{n_k})}}}(x)+(1-c),\ \mbox{ almost surely for }x> 0.
\end{align*}
By equating both the representations of the limit, we get that $F_{\frac{A_{n_k}}{\sqrt{m_{n_k}(\mathbf p_{n_k})}}}(x)$ converges to $F_{\mu_{sc}}(x)$ for all $x \neq 0$  and therefore
$\mu_{\frac{A_{n_k}}{\sqrt{m_{n_k}(\mathbf p_{n_k})}}} \stackrel{w}{\longrightarrow} \mu_{sc}$ almost surely. 

Now due to Lemma \ref{lem:subsequence_conv}, $\mu_{\frac{A_{n}}{\sqrt{m_{n}(\mathbf p_{n})}}} \stackrel{w}{\longrightarrow} \mu_{sc}$ in probability, and this complete the proof for unsigned adjacency matrix case. The proof for signed case is similar and we skip the details. 
\end{proof}

\section{Semicircle law for Upper model}\label{sec:upper_model}


Recall the upper model $Z_d(n,\mathbf{p})$ defined in \eqref{eqn:upper_model_defn}. Note that for $\sigma, \sigma^\prime \in K^{d-1}$ such that $\sigma \cup \sigma^\prime \in K^{d}$, $\sigma \cup \sigma^\prime \in Z_d(n,\mathbf{p})$  if and only if $\sigma \cup \sigma^\prime$ belongs the hypergraph $\mathcal{X}(n,\mathbf{p})$. Therefore, for each $\sigma, \sigma^\prime \in K^{d-1}$ such that $\tau= \sigma \cup \sigma^\prime \in K^d$, the distribution of the corresponding entry in the extended adjacency matrix, $\widehat{A}_{\sigma \sigma^\prime}$, has the distribution $\Ber{p_d}$, and is independent of entries $\widehat{A}_{\widetilde{\sigma}\widetilde{\sigma}^\prime}$ for which $\widetilde{\sigma} \cup \widetilde{\sigma}^\prime \neq \tau$. 
Hence, the extended unsigned adjacency matrix $\widehat{A}_n$ of the upper model is exactly the matrix $A_n(d)$ given in \eqref{eqn:adj_hadamard}, the adjacency matrix of the Linial-Meshulam complex with parameter $p_d$, indexed by $K^{d-1}$. For the multi-parameter upper model $Z_d(n,\mathbf{p})$, we first prove the following lemma.
\begin{lemma}\label{lemma:upper_model_dim}
    For fixed $d \geq 2$ and all choices of $(p_1,p_2,\ldots,p_d)$ such that $np_d \to \infty$, 
    \[
    \frac{|Z^{d-1}|}{\binom{n}{d}} \to 1\ \mbox{almost surely.}
    \]
\end{lemma}
\begin{proof}
We prove the cases $d \geq 3$ and $d=2$ separately. First, we look at the case $d \geq 3$. Note that 
    \[1 \geq\frac{|Z^{d-1}|}{\binom{n}{d}} \geq \frac{| \{\sigma \in Z^{d-1}: \sigma \subset \tau \text{ for some } \tau \in Z^d\}|}{\binom{n}{d}}, \]
    and we will prove that the fraction on the right converges to 1 almost surely.

    For convenience, we denote \[
    T_n:= |\{\sigma \in Z^{d-1}: \sigma \subset \tau \text{ for some } \tau \in Z^d\}|.
    \]
    We first note that  for all fixed $\sigma_0 \in K^{d-1}$,
    \[
    \mathbb{P}(\sigma_0  \subset \tau \text{ for some } \tau \in Z^d)=1-(1-p_d)^{n-d} \geq 1-e^{-np_d},
    \]
    where the last inequality follows since $(1-x) \leq e^{-x}$ for all real $x$. Since $np_d \to \infty$, it follows that for all $\sigma_0 \in K^{d-1}$,
     $\mathbb{P}(\sigma_0  \subset \tau \text{ for some } \tau \in Z^d) \to 1$ as $n \to \infty$. 
Now, consider the family of random variables $(X_{\sigma})_{\sigma \in K^{d-1}}$ given by 
\[
X_\sigma=\mathbf{1}(\sigma  \subset \tau \text{ for some } \tau \in Z^d)\]
and note that $T_n=\sum_{\sigma\in K^{d-1}} X_{\sigma}$. We calculate the variance of $T_n$.

    Observe that when $|\sigma \cap \sigma^\prime|< d-1$, the random variables  $X_{\sigma}$ and $X_{\sigma^\prime}$ are independent, and therefore $\Cov(X_{\sigma},X_{\sigma^\prime})=0$. Hence 
    \[
    \var(T_n)= \sum_{\sigma \in K^{d-1}}\var(X_\sigma)+ \sum_{\sigma,\sigma^\prime \in K^{d-1} :\atop |\sigma \cap \sigma^\prime|=d-1} \Cov(X_{\sigma},X_{\sigma^\prime}).
    \]
    Now, since  $|\Cov(X_{\sigma},X_{\sigma^\prime})|\leq 2$,  the first sum in the right hand side is of the order $O(n^d)$ and the second summation is of the order $O(n^{d+1})$. As a result, we get that $\var\left(\frac{T_n}{\binom{n}{d}}\right)=O(n^{-d+1})$, and consequently $\frac{T_n}{\binom{n}{d}}$ converges to 1 almost surely for all $d \geq 3$. This completes the proof for the case $d \geq 3$. 
    
    We prove the $d=2$ case by using more robust bounds on covariance and a coupling argument. For this, we note that $\var(X_\sigma) \leq 1$ for all $\sigma$ and for $\sigma \neq \sigma'$,
    \begin{align*}
     \Cov(X_\sigma,X_\sigma^\prime)&
    = \mathbb{P}(\sigma\subset \tau, \sigma'\subset \tau' \mbox{for some}\  \tau, \tau'\in Z^2 )-(\mathbb{P}(\sigma\subset \tau \text{ for some } \tau\in Z^2))^2\\
    &=[1-2(1-p_2)^{n-2}+(1-p_2)^{2n-5}]-[1-(1-p_2)^{n-2}]^2\\
    &=p_2(1-p_2)^{2n-5}.
    \end{align*}
    Therefore 
    \[
    \var\left(\frac{T_n}{\binom{n}{2}}\right) = \frac{1}{\binom{n}{2}^2} \sum_{\sigma \in K^{1}}\var(X_\sigma)+ \frac{1}{\binom{n}{2}^2}\sum_{\sigma,\sigma^\prime \in K^{1} :\atop |\sigma \cap \sigma^\prime|=1} \Cov(X_{\sigma},X_{\sigma^\prime}) \leq \frac{1}{\binom{n}{2}}+\frac{2(n-2)p_2}{\binom{n}{2}}.
    \]
    Choosing
    $p_2^{(n)}=\frac{1}{\sqrt{n}}$, we get that in this case $\sum_{n} \var\left(\frac{T_n}{\binom{n}{2}}\right) <\infty$ and therefore $ \frac{| Z^{1}(p_2^{(n)})|}{\binom{n}{2}} \to 1$ almost surely.
    
For an arbitrary choice of $p_2$, we argue as follows: Consider a sequence $p_2^{(n)}$ such that $np_2^{(n)} \to \infty$. We consider the following partition of this sequence: 
    \[
    A_1=\{ n \in \N: p_2^{(n)} \leq \frac{1}{ \sqrt{n}}\} \quad \text{and} \quad A_2=\{ n \in \N: p_2^{(n)} > \frac{1}{ \sqrt{n}}\}.
    \]
    Note that along the subsequence $A_1$, $\sum \var\left(\frac{T_n}{\binom{n}{2}}\right) <\infty$ and therefore $ \frac{| Z^{1}(\widetilde{p_2}^{(n)})|}{\binom{n}{2}} \to 1$ almost surely. And to prove the result for the subsequence in $A_2$, we consider $\widetilde{p_2}^{(n)}=\frac{1}{\sqrt{n}}$ and couple $Z^2(p_2^{(n)})$ and $Z^2(\widetilde{p_2}^{(n)})$. Note that the $2$-cells of the random simplicial complexes are constructed from the collection of random variables $(X_{\tau})_{\tau \in K^2}$.  Consider a sequence of independent uniform distributions $U_\tau \sim U[0,1]$, $\tau \in K^2$ and let 
    \[
    X_{\tau}^{(n)}=\mathbf{1}(U_\tau \leq p_2^{(n)}) \quad \text{and} \quad \widetilde{X}_{\tau}^{(n)}=\mathbf{1}(U_\tau \leq \widetilde{p_2}^{(n)}).
    \]
     Note that under this coupling 
    \begin{multline*}
         | \{\sigma \in Z^{1}(\widetilde{p_2}^{(n)}): \sigma \subset \tau \text{ for some } \tau \in Z^2(\widetilde{p_2}^{(n)})\}| \\\leq  |\{\sigma \in Z^{1}(p_2^{(n)}): \sigma \subset \tau \text{ for some } \tau \in Z^2(p_2^{(n)})\}|.
    \end{multline*}
    Now, since $$\frac{| \{\sigma \in Z^{1}(\widetilde{p_2}^{(n)}): \sigma \subset \tau \text{ for some } \tau \in Z^2(\widetilde{p_2}^{(n)})\}| }{\binom{n}{2}} \to 1 \ \mbox{ almost surely,}$$ it follows that 
    $$ \frac{|\sigma \in Z^{1}(p_2^{(n)}): \sigma \subset \tau \text{ for some } \tau \in Z^2(p_2^{(n)})\}|}{\binom{n}{2}} \to 1 \ \mbox{ almost surely.}$$ This completes the proof.
\end{proof}
 
Now we prove Theorem \ref{thm:semicircle_upper}.
\begin{proof}[Proof of Theorem \ref{thm:semicircle_upper}]
    Recall the observation that the extended unsigned adjacency matrix $\widehat{A}_n$ of $Z_d(n,\mathbf p)$ is the same as the unsigned adjacency matrix of Linial-Meshulam complex. Note from \eqref{eqn:block_form} that $$\mu_{\widehat{A}_n}(\{0\})= \frac{\dim(A_n)}{\binom{n}{d}}\mu_{A_n}(\{0\})+1-\frac{\dim(A_n)}{\binom{n}{d}}.$$ Note that $np_d(1-p_d)\to\infty$ implies that $np_d\to\infty$. Therefore, by Lemma \ref{lemma:upper_model_dim} we have that $\frac{\dim(A_n)}{\binom{n}{d}} \to 1$ in probability, as $n \to \infty$, and hence $\mu_{A_n}(\{0\})$ converges to 0 in probability. 
Using Result \ref{res:LSD_linialMeshulam}, we get that, for $x<0$,
\begin{align*}
\lim_{n \to \infty }F_{\frac{A_n}{\sqrt{np_d(1-p_d)}}}(x)&=\lim_{n \to \infty }\frac{\binom{n}{d}}{\dim(A_n)}F_{\frac{\widehat{A}_n}{\sqrt{np_d(1-p_d)}}}(x)\\
&= F_{\mu_{sc}}(x)  
\end{align*}
and for $x\geq 0$,
\begin{align*}
\lim_{n \to \infty }F_{\frac{A_n}{\sqrt{np_d(1-p_d)}}}(x)&=\lim_{n \to \infty } \Big[\frac{\binom{n}{d}}{\dim(A_n)}F_{\frac{\widehat{A}_n}{\sqrt{np_d(1-p_d)}}}(x)+1-\frac{\binom{n}{d}}{\dim(A_n)}\Big]  \\&= \lim_{n \to \infty }\frac{\binom{n}{d}}{\dim(A_n)}F_{\frac{\widehat{A}_n}{\sqrt{np_d(1-p_d)}}}(x)\\
&= F_{\mu_{sc}}(x),
\end{align*}
in probability, where $F_{\mu_{sc}}$ denotes the cumulative distribution function of the  semicircle law.  Proof for the signed case is similar and this completes the proof.
\end{proof}

\section{Counting $(d-1)$-cells}\label{sec:count_simplices}
In this section, we prove Lemma \ref{lem:num_simplices} and Theorem \ref{thm:maximal_faces}. The proof ideas of these theorems are based on the techniques in \cite{Fowler_MultiRSC_homology_2019}.

\begin{proof}[Proof of Lemma \ref{lem:num_simplices}]
Let $f_{d-1}=\sum_{\sigma \in K^{d-1}}\mathbf{1}_Y(\sigma)$ denote the number of $(d-1)$-cells in $Y=\RSC$. 
We begin by calculating the expected number of $(d-1)$-cells. We have
\[
\mathbb{E}[f_{d-1}] = \binom{n}{d} \prod_{i=1}^{d-1} p_i^{\binom{d}{i+1}}, 
\]
and since $(\RSC)_n$ belongs to $\domc$, we get
\begin{equation}\label{eqn:conv_f_d-1}
  \lim_{n \to \infty}\frac{\mathbb{E}[f_{d-1}]}{{n \choose d}} = c.  
\end{equation}
Next, we calculate $\mathbb{E}[f_{d-1}^2]$. For $\sigma \in K^{d-1}$, let $E_\sigma$ be the event that $\sigma \in \RSC$. Then
\[
\mathbb{E}[f_{d-1}^2] = \sum_{\sigma,\sigma^\prime \in K^{d-1}} \mathbb{P}[E_\sigma \cap E_\sigma^\prime]
= \binom{n}{d} \sum_{\sigma^\prime \in K^{d-1}} \mathbb{P}[E_\sigma \cap E_{\sigma^\prime}], 
\]
where the second equality follows from the homogeneity of the model.

 We proceed by grouping together all $\sigma^\prime \in K^{d-1}$ for which the number of vertices in $\sigma \cap \sigma^\prime$ is the same. Note that if $\sigma$ and $\sigma^\prime$ have $m$ vertices in common, then the number of ways of choosing $\sigma^\prime$ is ${d \choose m }{n-d \choose d-m}$ and the total number of $i$-cells in $\sigma \cap \sigma^\prime$ is ${m \choose i+1}$.
 We therefore have
\begin{align*}
 \mathbb{E}[f_{d-1}^2] &=
\binom{n}{d} \sum_{m=0}^d \binom{d}{m} \binom{n-d}{d-m} \prod_{i=1}^{d-1} p_i^{2\binom{d}{i+1} - \binom{m}{i+1}} \\
&=\sum_{m=0}^d \binom{n}{d} \binom{d}{m} \binom{n-d}{d-m}\prod_{i=1}^{d-1} p_i^{2\binom{d}{i+1}}    \prod_{i=1}^{m-1} p_i^{-\binom{m}{i+1}}.   
\end{align*}

Note that for $m=0$, the term inside the summation can be bounded as
\[
\binom{n}{d} \prod_{i=1}^{d-1} p_i^{2\binom{d}{i+1}}   \binom{n-d}{d} \leq \left(\binom{n}{d}\prod_{i=1}^{d-1} p_i^{\binom{d}{i+1}} \right)^2=\E[f_{d-1}]^2.
\]
Therefore we get that 
\[
\mathbb{E}[f_{d-1}^2] -\mathbb{E}[f_{d-1}]^2 \leq \binom{n}{d} \prod_{i=1}^{d-1} p_i^{2\binom{d}{i+1}}
\left( \sum_{m=1}^d \binom{d}{m} \binom{n-d}{d-m} \prod_{i=1}^{m-1} p_i^{-\binom{m}{i+1}} \right).
\]

We consider the terms corresponding to each $m$ separately, and we have for $1 \leq m \leq d$,
\[
\binom{n}{d}\binom{d}{m} \binom{n-d}{d-m} \prod_{i=1}^{m-1} p_i^{-\binom{m}{i+1}}=O\left( n^d \times n^{d-m} \right)= O(n^{2d-m}),
\]
where the first equality follows since $\liminf p_i>0$ for all $i$. Therefore, we get that $\mathbb{E}[f_{d-1}^2] -\mathbb{E}[f_{d-1}]^2$ is of the order $O(n^{2d-1})$.
As a result, we have
\[
\operatorname{Var}\left[\frac{f_{d-1}}{{n \choose d}}\right] = \frac{\mathbb{E}[f_{d-1}^2] - \mathbb{E}[f_{d-1}]^2}{{n \choose d}^2}
= O(n^{-1}),
\]
and from \eqref{eqn:conv_f_d-1}, we get that $\frac{f_{d-1}}{{n \choose d}}$ converges  to $c$ in probability.
\end{proof}
\begin{proof}[Proof of Theorem \ref{thm:maximal_faces}]
We first prove (i). Recall that $N_{d-1}(Y)$ denotes the number of maximal $(d-1)$-cells of $\RSC$, and note that
\begin{align*}
   \mathbb{E}[N_{d-1}(Y)] 
&= \binom{n}{d} \left( \prod_{i=1}^{d-1} p_i^{\binom{d}{i+1}} \right) \left( 1 - \prod_{i=1}^d p_i^{\binom{d}{i}} \right)^{n-d},
\end{align*}
and therefore
\begin{align*}
  \limsup_{n \to \infty} \mathbb{E}[N_{d-1}(Y)]    &\leq \limsup_{n \to \infty}  \binom{n}{d}
\exp\!\left( -(n-d)\prod_{i=1}^d p_i^{\binom{d}{i}} \right),
\end{align*}
using the fact $(1-x) \leq e^{-x}$ for all $x$, and $p_i \leq 1$ for all $i$.
Then there exists some constant $D > 0$ independent of $n$,
\[
\mathbb{E}[N_{d-1}(Y)] \leq D
\exp\left(d\log n - np_d \prod_{i=1}^{d-1} p_i^{d \choose i} \right).
\]
By hypothesis
\[
\min_{i=1,\ldots, d-1}\liminf p_i^{d \choose i}>0 \text{ and } p_d \gg \frac{\log n}{n},
\]
Therefore we have that for sufficiently large $n$, we have $np_d\prod_{i=1}^{d-1} p_i^{d \choose i}>(d+2)\log n$ and therefore
\[
 D
\exp\left(d\log n - np_d \prod_{i=1}^{d-1} p_i^{d \choose i} \right) \leq D \exp \left(- n \times \frac{2\log n }{n}\right)= O(n^{-2}).
\]
 Now, using Markov’s inequality, it follows that
\[
\mathbb{P}[N_{d-1}(Y) \geq 1] \leq \mathbb{E}[N_{d-1}(Y)] = O(n^{-2}),
\]
and part (i) of the theorem follows as a consequence of Borel-Cantelli lemma. Now, part (ii) of the theorem follows from observing that $N_{d-1}(Z) = N_{d-1}(Y)$.
\end{proof}
\textbf{Acknowledgment.} Kiran Kumar acknowledges the helpful discussions with Manjunath Krishnapur during the initial stages of the work. In particular, his suggestion to consider the extended adjacency matrix as a Hadamard product. The research of Kartick Adhikari was partially supported by the Inspire Faculty Fellowship: DST/INSPIRE/04/2020/000579
\bibliography{ref}
\bibliographystyle{amsplain}
\end{document}